\documentclass[final,3p,times]{elsarticle}
\usepackage{graphicx}
\usepackage{subfig,comment}
\usepackage{mathptmx}      
\usepackage{xcolor}
\usepackage{caption}

\usepackage{amsfonts,amsthm}
\usepackage{amsmath,amssymb,amssymb,amscd,,subeqnarray}

\newtheorem{remark}{Remark}

\usepackage{graphicx}%
\usepackage{boldfonts}
\usepackage{siunitx} 

\usepackage{stmaryrd}
\usepackage{listings}
\usepackage{boldfonts}
\usepackage{symbols}
\usepackage{multirow}
\usepackage{cases}
\usepackage{tikz}
\usetikzlibrary{shapes.geometric}
\usetikzlibrary{shapes,arrows}

\tikzstyle{block} = [rectangle, draw, text centered, rounded corners, 
                     text width=7em,    
                     minimum height=3em, node distance=3.3cm]
\tikzstyle{line} = [draw, -latex']
\tikzstyle{cloud} = [draw, ellipse, node distance=2.5cm, minimum height=3em]

\newcommand{\vertiii}[1]{{\left\vert\kern-0.25ex\left\vert\kern-0.25ex\left\vert #1 
    \right\vert\kern-0.25ex\right\vert\kern-0.25ex\right\vert}}
\newcommand{\bdf}[2]{{\textup{\textsf{BDF}}_{#2}({#1})}}

\journal{Journal of Computational Physics}

\begin{document}

\begin{frontmatter}



\title{Enriched Galerkin methods for two-phase flow in porous media with capillary pressure}


\author[shlee]{Sanghyun Lee\corref{cor1}}
\ead{lee@math.fsu.edu}

\author[mfw]{Mary F. Wheeler}
\ead{mfw@ices.utexas.edu}

\cortext[cor1]{Corresponding author}


\address[shlee]{
Department of Mathematics, 
Florida State University, 1017 Academic Way, Tallahassee, FL 32306-4510, United States}

\address[mfw]{The Center for Subsurface Modeling, The Institute for Computational Engineering and Sciences,The University of Texas at Austin, Austin, TX 78712, United States}

\begin{abstract}
In this paper, we propose {an} enriched Galerkin (EG) approximation
for {a two-phase pressure saturation system with capillary pressure}  in heterogeneous porous media.
The EG methods  are locally conservative, have {fewer} degrees of freedom compared to  discontinuous Galerkin (DG),
and have an efficient pressure solver.
To avoid non-physical oscillations, 
an entropy viscosity stabilization method is employed for {high} 
order saturation {approximations}.  
Entropy residuals {are} applied for dynamic mesh adaptivity to reduce the computational cost for larger computational domains.  
The iterative {and sequential} IMplicit Pressure and Explicit Saturation (IMPES) algorithms are {treated in time}.
Numerical examples with different relative permeabilities and capillary pressures are {included} to verify and to {demonstrate} the capabilities of  EG.
\end{abstract}

\begin{keyword}
Enriched Galerkin finite element methods 
\sep Two-phase flow 
\sep Capillary pressure 
\sep Porous media 
\sep Entropy viscosity 
\sep Dynamic mesh adaptivity
\end{keyword}
\end{frontmatter}

\section{Introduction}\label{sec:intro}
We consider a two-phase flow system in porous media which has been widely employed in petroleum reservoir modeling and environmental engineering for the 
past several decades  \cite{aziz1979petroleum,chavent1986mathematical,Douglas:1979hg,morel1973two,peaceman2000fundamentals,slattery1970two}.
The conventional two-phase flow system is formulated by coupling Darcy's law  for multiphase flow with the saturation transport equation \cite{kueper1991two,whitaker1986flow}. 

An incomplete list of numerical approximations
such as finite difference, mixed finite elements, and finite volume methods 
\cite{
arbogast1992existence,
arbogast2002implementation,
aziz1979petroleum,
coats1982reservoir,
coats2000note,
Douglas:1979hg,
EpstheynRiviere2007-apnum,
EphsteynRiviere2009,
peaceman2000fundamentals,
radu2015robust,
riaz2006numerical,
Schmid2013416,
Zhang201339} 
have been successfully utilized in  multiphase flow reservoir simulators.
Recent interest has centered on multiscale  extensions to finite element methods  \cite{
arbogast2000numerical,
efendiev2006accurate,
efendiev2003accurate,
ganis2014global,
Hajibeygi20095129,
jenny2003multi,
lee2008multiscale,Peszynska:2002bc}. 
In all of these works, it was observed that local conservation was required for accurately solving the saturation transport equations \cite{KAASSCHIETER1995277,ScovazziWheMikLee16}. 
However, only several of these references considered capillary pressure effects for two-phase flow systems  
\cite{Arbogast:2013eg,
Bastian:2014fn,
ELAMIN20161344,
Ern:2010cn,
Hoteit:2008dt,
kou2010new, 
riaz2004linear, yang2017nonlinearly}. 
For many problems such as CO$_2$ sequestration, the latter is crucial for realistic heterogeneous media. 

In this paper, we  focus on extensions of  enriched Galerkin approximations (EG) {to} two-phase flow in porous media with capillary pressure. 
Our objective is to demonstrate that high order spatial approximations for saturations can be computed efficiently using EG.
EG provides locally and globally conservative fluxes and preserves local mass balance for transport \cite{LeeLeeWhi15, LeeLeeWhi17a,LeeWhe_2017}.
EG is constructed by enriching the
conforming continuous Galerkin finite element method (CG) with piecewise constant functions \cite{BecBurHansLar2003,sunliu2009},
with the same bilinear forms as the interior penalty DG schemes.
However, EG has substantially fewer degrees of freedom in comparison with DG 
and a fast effective high order solver for pressure whose cost is roughly that of CG  \cite{LeeLeeWhi15}.
EG has been successfully employed to realistic multiscale and multi-physics applications \cite{LeeWhe_2017,LeeMiWheWi17,LeeMiWheWi16}.
An additional advantage of EG is that 
only those subdomains that require local conservation need be enriched with a treatment of high order non-matching grids. 

Local conservation of the flux is crucial for flow and 
saturation stabilization is critical for avoiding overshooting, undershooting, and spurious oscillations \cite{kou2013convergence}. 
Our high order EG transport system is coupled with an entropy viscosity residual stabilization method introduced in \cite{guermond2011entropy} to 
avoid spurious oscillations near the interface of {saturation fronts}.
Instead of using limiters and non-oscillatory reconstructions,  
this method adds nonlinear dissipation to the numerical discretization \cite{Guermond2015,Guermond2017448,Guermond2011}. 
The numerical diffusion is constructed by the local residual of an entropy residual. 
Moreover, the entropy residual is employed for dynamic adaptive mesh refinement to capture the moving interface between the immiscible fluids 
\cite{jenny2005adaptive,KlieberRiviere2006}.  
It is shown in \cite{andrews1998posteriori,puppo2004numerical} that the entropy residual can be used as an a posteriori error indicator.

To take advantage of high order in space, each time derivative in the flow and transport system is discretized by second order backward difference formula (BDF2) and extrapolations are employed. 
For the coupling solution algorithm, 
{a sequential time-stepping} scheme (IMPES)  is applied for efficient computation \cite{ewing1983simulation}.
First, we solve the pressure equation implicitly 
assuming saturation values are obtained by extrapolation in time and 
 the transport equation is solved explicitly  
 \cite{chen2004improved,fagin1966new, kou2010new,kou2010iterative, lu2009iterative, young1983generalized}. 
In addition, we employ H(div) flux reconstruction to the incompressible flow to enhance the performance as applied for DG in \cite{bastian2003superconvergence,ern2007accurate,Li:2015fl}.

\section{Mathematical Model}
\label{sec:Model}
In this section, a mathematical model for the slightly compressible two-phase Darcy flow and saturation system in a heterogeneous media is presented. 
Let $\Omega\subset \Reals{d}$ be a bounded polygon (for $d=2$) or polyhedron (for $d=3$) with Lipschitz boundary $\partial \Omega$, and $(0,\mathbb{T}]$  the computational time interval with $\mathbb{T} > 0$.
{The mass conservation equation for saturation equation is defined by}
\begin{equation}
\dfrac{\partial}{\partial t} (\phi \rho_i s_i ) 
 + \nabla \cdot ( \rho_i \bu_i ) = \rho_i f_i, \; i \in \{w,n\},
\label{eqn:main_s}
\end{equation}
where $\phi$ is the porosity of the porous media, $\rho_i$ is the density, $s_i : \Omega \times (0,\mathbb{T}] \rightarrow \mathbb{R}$ is the saturation, 
and $i \in \{w,n\}$ indicates wetting$(w)$ or non-wetting$(n)$ phases, respectively. 
{Here, $f_i :=  \tilde{s}_i q_i$,
where $\tilde{s}_i, q_i$ are the saturation injection/production term and flow  injection/production, respectively.
If $q_i>0$, $\tilde{s}_i$ is the injected saturation of the fluid and 
if $q_i<0$, $\tilde{s}_i$ is the produced saturation.}
Here $\bu_i : \Omega \times (0,\mathbb{T}] \rightarrow \mathbb{R}^d$ is the Darcy velocity for each phase i, given by 
\begin{equation}
\bu_i := -  \bK \dfrac{k_i}{\mu_i} \left( \nabla p_i - \rho_i \bg \right),
\label{eqn:darcy}
\end{equation}
in which $k_i$ is the relative permeability,
$\bK:=\bK(\bx)$ is the absolute permeability tensor of the porous media,
$\mu_i$ is the viscosity, 
$p_i : \Omega \times (0,\mathbb{T}] \rightarrow \mathbb{R}$ is the pressure for each phase, and $\bg$ is the gravity acceleration.
Relative permeability is a given function of saturation which is defined as 
\begin{equation}
k_i := k_i(s_w).
\end{equation}
Here we define the capillary pressure, 
\begin{equation}
p_c := p_c(s_w) = p_n - p_w,
\label{eqn:capillary_pressure}
\end{equation}
which is the pressure difference between the wetting and non-wetting phase \cite{ChenHuan2006}. 
Since, we assume that all pores are filled with fluid, we have
\begin{equation}
s_w + s_n = 1 
\ 
\text{ and } 
\ 
\tilde{s}_w + \tilde{s}_n = 1.
\end{equation}

To derive a pressure equation, we sum the saturation equations \eqref{eqn:main_s} to get 
\begin{equation}
\phi \dfrac{\partial}{\partial t} \left(\rho_w s_w + 
\rho_n s_n \right) 
+ \nabla \cdot ({\rho_w}\bu_w +{\rho_n}\bu_n)  = \rho_w f_w + \rho_n f_n,
\label{eqn:slightly_comp}
\end{equation}
where we consider a slightly compressible fluid satisfying 
\begin{equation}
\rho_i(p_i) 
\approx  \rho_i^0 \exp^{c_i^F(p_i-p_i^0)}
\approx \rho_i^0(1+c_i^F(p_i-p_i^0)), 
\end{equation}
with a small compressibility coefficient, $c_i^F \ll 1$. 
Here we assume the reference pressure $p_i^0$ is zero, and porosity $\phi$ and reference density $\rho_i^0$ are constants. 
Thus, we can rewrite \eqref{eqn:slightly_comp} and obtain  
\begin{equation}
\phi \dfrac{\partial}{\partial t} \left(c_w^F \rho_w^0 p_w s_w  + c_n^F \rho_n^0 p_n s_n \right) 
+ \nabla \cdot ({\rho_w}\bu_w +{\rho_n}\bu_n)  = \rho_w f_w + \rho_n f_n.
\label{eqn:slightly_comp_flow}
\end{equation}
For the incompressible case, we set $c_i^F=0$ and have 
\begin{equation}
\nabla \cdot ({\rho_w}\bu_w +{\rho_n}\bu_n)  =  \rho_w f_w + \rho_n f_n.
\label{eqn:incomp}
\end{equation}

\subsection{Choice of primary variables}
Throughout the paper, we set the wetting phase pressure $p_w$  
and saturation $s_w$ as the primary variables. 
Different choices and effects are illustrated in \cite{Arbogast:2013eg}. 
We rewrite the incompressible flow equation 
by combining the relations  \eqref{eqn:darcy}, \eqref{eqn:capillary_pressure},  \eqref{eqn:incomp}, and continuity of phase fluxes to obtain 
\begin{equation}
-\nabla \cdot (\bK\lambda_t ( \nabla p_w
 - \rho_w \bg )
 + \lambda_n (\bK \nabla p_c + (\rho_w - \rho_n)\bg )) 
 = (\rho f)_t,
\label{eqn:pressure_00}
\end{equation}
which is equivalent with 
\begin{equation}
-\nabla \cdot (\bK (\lambda_t  \nabla p_w
 - (\rho \lambda)_t \bg )
 + \bK \lambda_n  \nabla p_c) = (\rho f)_t,
\label{eqn:pressure}
\end{equation}
where 
\begin{align}
\lambda_i &:= \lambda_i(s_w) =  \rho_i\dfrac{k_i(s_w)}{\mu_i} , \; \text{phase mobility}\\
\lambda_t &:= \lambda_t(s_w) = \lambda_w(s_w) + \lambda_n(s_w), \; \text{total mobility}\\
(\rho \lambda)_t &:= (\rho \lambda(s_w))_t  = \rho_w \lambda_w(s_w) + \rho_n \lambda_n(s_w), \\
(\rho f)_t  &:=\rho_w f_w + \rho_n f_n. 
\end{align}
For the slightly compressible flow equations, we get the pressure equation
\begin{multline}
\label{eqn:slight_comp}
\phi \dfrac{\partial}{\partial t} 
\left( c_w^F \rho_w^0  s_w  p_w  
+  c_n^F \rho_n^0  (1-s_w)  p_w  + c_n^F \rho_n^0  (1-s_w) p_c \right) 
-\nabla \cdot (\bK (\lambda_t  \nabla p_w
 - (\rho^0 \lambda)_t \bg )
 + \bK  \lambda_n \nabla p_c) = (\rho^0 f)_t,
\end{multline}
where 
\begin{align}
(\rho^0 \lambda)_t &:= \rho^0_w \lambda_w + \rho^0_n \lambda_n, \\
(\rho^0 f)_t  &:=\rho^0_w f_w + \rho^0_n f_n. 
\end{align}
For the saturation equation,  we solve 
\begin{equation}
\dfrac{\partial}{\partial t} (\phi \rho^0_w s_w ) 
 + \nabla \cdot ({\rho^0_w} \bu_w ) = \rho^0_w f_w,
 \label{eqn:transport}
\end{equation}
and $s_w +s_n =1$.

The boundary of $\Omega$ is decomposed into three disjoint sets
$\Gamma_\textsf{in}$, $\Gamma_\textsf{out}$ and $\Gamma_N$ so that 
$\overline{\partial \Omega} = 
\overline{\Gamma}_\textsf{in}
\cup
\overline{\Gamma}_\textsf{out}
\cup \overline{\Gamma}_N$ 
For the flow problem,  we impose
\begin{align}
p_w (\text{ or } p_n)  = p_\textsf{in}  
\mbox{ on } \Gamma_\textsf{in} 
&\times (0,\mathbb{T}], \quad \\ 
p_w (\text{ or } p_n)  = p_\textsf{out}  
\mbox{ on } \Gamma_\textsf{out} 
&\times (0,\mathbb{T}], \quad \\ 
(\bu_w + \bu_n) \cdot \bn = \bu_{N}  
\mbox{ on } \Gamma_N 
&\times (0,\mathbb{T}], 
\end{align}
where $p_\textsf{in} \in L^2(\Gamma_\textsf{in})$,  
$p_\textsf{out} \in L^2(\Gamma_\textsf{out})$ and $\bu_N \in L^2(\Gamma_N)$ are the each Dirichlet and Neumann boundary conditions, respectively. Thus we define ${\Gamma}_D := \Gamma_\textsf{in} \cup \Gamma_\textsf{out}$.
Here inflow and outflow boundaries are defined as 
$$
\Gamma_{\rm in} := \{ \bx \in \partial \Omega : \bu_w\cdot\bn < 0\} \ \text{ and }
\Gamma_{\rm out} := \{ \bx \in \partial \Omega : \bu_w\cdot\bn > 0\}. 
$$
For the saturation system, we impose
\begin{equation}
s_w (\text{ or } s_n) = 
s_\textsf{in}, 
 \mbox{ on } \Gamma_\textsf{in} \times (0,\mathbb{T}]
\label{main_bd_c}
\end{equation}
where $s_\textsf{in}$ is a given boundary value for saturation. 
Finally, the above systems are supplemented by initial conditions 
$$
s_w(\bx,0) = s_w^0(\bx), \mbox{ and }
p_w(\bx,0) = p_w^0(\bx), \quad \forall \bx \in \Omega.
$$

\section{Numerical Method}
\label{sec:discret}
Let $\mathcal{T}_h$ be the shape-regular (in the sense of Ciarlet)  triangulation by a family of partitions of $\O$ into $d$-simplices $\K$ (triangles/squares in $d=2$ or tetrahedra/cubes in $d=3$). 
We denote by $h_{\K}$ the diameter of $\K$ and we set $h=\max_{\K \in \Th} h_{\K}$.  
Also we denote by $\Eh$ the set of all edges and by $\Eho$ and $\Ehb$ the collection of all interior and boundary edges, respectively. 
In the following notation, we assume edges for two dimension but the results hold analogously for faces in three dimensional case.
For the flow problem, the boundary edges $\Ehb$ can be further decomposed into $\Ehb = \mathcal{E}_h^{D,\partial} \cup \mathcal{E}_{h}^{N,\partial}$, where $\mathcal{E}_h^{D,\partial}$ is the collection of edges where the Dirichlet boundary condition is imposed (i.e $\mathcal{E}_h^{D,\partial} := 
\mathcal{E}_h^{\textsf{in},\partial} \cup 
\mathcal{E}_{h}^{\textsf{out},\partial}$), 
while $\mathcal{E}_h^{N,\partial}$ is the collection of edges where the Neumann boundary condition is imposed. 
In addition, we let  $\mathcal{E}_h^{1} := \Eho \cup \mathcal{E}_h^{D,\partial}$ and $\mathcal{E}_h^{2} := \Eho \cup \mathcal{E}_h^{N,\partial}$. 
For the transport problem,  
the boundary edges $\Ehb$ decompose into $\Ehb = \mathcal{E}_h^{\text{in}} \cup \mathcal{E}_{h}^{\text{out}}$, where $\mathcal{E}_h^{\text{in}}$ is the collection of edges where the inflow boundary condition is imposed, while $\mathcal{E}_h^{\text{out}}$ is the collection of edges where the outflow boundary condition is imposed.

The space $H^{s}(\Th)$ $(s\in \Reals{})$ is the set of element-wise $H^{s}$ functions on $\mathcal{T}_h$, and $L^{2}(\Eh)$ refers to the set of functions whose traces on the elements of $\Eh$ are square integrable. Let $\mathbb{Q}_l(\K)$ denote the space of polynomials of partial degree at most $l$. 
Regarding the time discretization, given an integer $N \geq 2$, we define a partition of the time interval 
$0 =: t^0 < t^1 < \cdots < t^N:= \mathbb{T}$ and denote $\Delta t := t^k - t^{k-1}$ for the uniform time step.
Throughout the paper, we use the standard notation for Sobolev spaces and their norms. For example, let $E \subseteq \Omega$, then $\|\cdot\|_{1,E}$ and $|\cdot|_{1,E}$ denote the $H^1(E)$ norm and seminorm, respectively. 
For simplicity, we eliminate the subscripts on the norms if $E = \Omega$.
For any vector space $\bX$, $\bX^d$ will denote the vector space of size d, whose components belong to $\bX$ and $\bX^{d\times d}$ will denote the $d \times d$ matrix whose components belong to $\bX$.

\begin{figure}[!h]
\centering
\begin{tikzpicture}
\node[] at (1,0.3) {${\K}_1$};
\node[] at (3,0.3) {${\K}_2$};

\draw (0,0) -- (2,0); 
\draw[dashed] (2,0) -- (2,2); 
\draw (2,2) -- (0,2); 
\draw (0,2) -- (0,0); 

\draw (2,0) -- (4,0); 
\draw (4,0) -- (4,2); 
\draw (4,2) -- (2,2); 
\draw[solid] (0,0) circle (0.05);
\draw[solid] (2.,0) circle (0.05);
\draw[solid] (2.,2) circle (0.05);
\draw[solid] (0,2) circle (0.05);
\draw[solid] (4,2) circle (0.05);
\draw[solid] (4,0) circle (0.05);
\node[] at (1,1) {\scriptsize $\triangle$} ; 
\node[] at (3,1) {\scriptsize $\triangle$} ;
\end{tikzpicture}
\caption{A sketch of the degrees of freedom for 
enriched Galerkin in 
a two-dimensional Cartesian grid ($\mathbb{Q}$) with $l=1$. 
Four circles ($\circ$) are the degrees of freedom for continous Galerkin $(M^l(\mathcal{T}_h))$ and ($\triangle$) is the discontinuous constant $(M^0(\mathcal{T}_h))$. }
\label{fig:EG}
\end{figure}
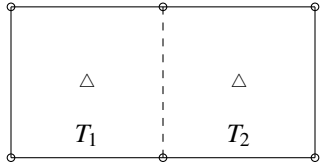

We  introduce the space of piecewise discontinuous polynomials of degree $l$ as
\begin{equation}
M^l(\mathcal{T}_h) := \left \{ \psi \in L^2(\Omega) | \ \psi_{|_{\K}} \in \mathbb{Q}_l(\K), \ \forall \K \in \mathcal{T}_h \right \}, 
\end{equation}
and let $M_0^l(\mathcal{T}_h)$ be the subspace of $M^l(\mathcal{T}_h)$ consisting of continuous piecewise polynomials;
\begin{equation*}
M_0^l(\mathcal{T}_h) = M^l(\mathcal{T}_h) \cap \mathbb{C}_0(\Omega). 
\end{equation*}
The enriched Galerkin finite element space, denoted by $V_{h,l}^{\textsf{EG}}$ is defined as
\begin{equation}
V_{h,l}^{\textsf{EG}}(\mathcal{T}_h)  := M^l_0(\mathcal{T}_h) + M^0(\mathcal{T}_h),
\end{equation}
where $l \geq 1$, also see \cite{BecBurHansLar2003,LeeLeeWhi15,LeeLeeWhi17a,LeeWhe_2017, sunliu2009} for more details.
\begin{remark}
We remark that the degrees of freedom for $V_{h,l}^{\textsf{EG}}(\mathcal{T}_h) $ when $l = 1$, is approximately one half and {one fourth} the degrees of freedom of the linear DG space, in two and three space dimensions, respectively. See Figure \ref{fig:EG}.
\end{remark}
We define the coefficient $\bkappa_T$ by 
\begin{equation}
\bkappa_T := \bkappa |_T, \quad \forall T \in \mathcal{T}_h.
\end{equation}
For any $e \in \Eho$, let $\K^{+}$ and $\K^{-}$ be two neighboring elements such that  $e = \partial \K^{+}\cap \partial \K^{-}$.
We denote by $h_{e}$ the length of the edge $e$. 
Let $\n^{+}$ and $\n^{-}$ be the outward normal unit vectors to  $\partial T^+$ and $\partial T^-$, respectively ($\n^{\pm} :=\n_{|\K^{\pm}}$). 
For any given function $\xi$ and vector function $\bf{\xi}$, defined on the triangulation $\mathcal{T}_h$, we denote $\xi^{\pm}$ and $\bxi^{\pm}$ by the restrictions of $\xi$ and $\bxi$ to $T^\pm$, respectively. 
We define the  average $\av{\cdot}$ as follows: for $\zeta \in L^2(\mathcal{T}_h)$ and $\taub \in L^2(\mathcal{T}_h)^d$,
\begin{equation}\label{av-w}
\av{\zeta} := \frac{1}{2} (\zeta^+ + \zeta^- )
\quad \mbox{ and } \quad 
\av{\taub} := \frac{1}{2} (\taub^+ +   \taub^-) \quad \mbox{on } e\in
\Eho.
\end{equation}
On the other hand, for $e \in \Ehb$, we set $\av{\zeta} :=   \zeta$ and $\av{\taub} :=  \taub$. 
The jump across the interior edge will be defined as usual: 
\begin{align*}
\jump{\zeta} = \zeta^+\n^++\zeta^-\n^- \quad \mbox{ and } \quad \jtau = \taub^+\cdot\n^+ + \taub^-\cdot\n^- \quad \mbox{on } e\in \Eho. 
\end{align*}
For inner products, we use the notations: 
\begin{align*}
&(v,w)_{\Th}:=\dyle\sum_{\K \in \Th} \int_{\K} v\, w dx, \quad \forall\,\, v ,w \in L^{2} (\mathcal{T}_h), \\
&\langle v, w\rangle_{\Eh}:=\dyle\sum_{e\in \Eh} \int_{e} v\, w \,d\gamma, \quad \forall\, v, w \in L^{2}(\Eh).
\end{align*}
For example, a function in 
$\psi_{\textsf{EG}} \in V_{h,l}^{\textsf{EG}}(\mathcal{T}_h)$ can be decomposed into 
$\psi_{\textsf{EG}} = \psi_{\textsf{CG}}+  \psi_{\textsf{DG}}$,
where 
$\psi_{\textsf{CG}} \in M^l_0(\mathcal{T}_h)$  
and
$\psi_{\textsf{DG}} \in M^0(\mathcal{T}_h)$. 
Thus the inner product $(\psi_{\textsf{EG}} ,\psi_{\textsf{EG}} ) = 
(\psi_{\textsf{CG}} ,\psi_{\textsf{CG}} ) + (\psi_{\textsf{CG}} ,\psi_{\textsf{DG}} ) + (\psi_{\textsf{DG}} ,\psi_{\textsf{CG}} ) + (\psi_{\textsf{DG}} ,\psi_{\textsf{DG}} )$ creates a matrix as 
$$
\begin{pmatrix} 
\psi_{\textsf{CG}} \psi_{\textsf{CG}} & \psi_{\textsf{CG}} \psi_{\textsf{DG}} \\
\psi_{\textsf{DG}} \psi_{\textsf{CG}} & \psi_{\textsf{DG}} \psi_{\textsf{DG}}
\end{pmatrix}.
$$
Finally, we introduce the interpolation operator $\Pi_h$ for the space $V_{h,l}^{\textsf{EG}}$ as
\begin{equation}
\Pi_h v = \Pi_0^l v + Q^0 ( v - \Pi_0^l v), 
\label{interp}
\end{equation}
where $\Pi_0^l$ is a {continuous} interpolation operator onto the space $M_0^l(\mathcal{T}_h)$, {and} $Q^0$ is the $L^2$ projection onto the space $M^0(\mathcal{T}_h)$.  See \cite{LeeLeeWhi15} for more details.

\subsection{Temporal Approximation}
The time discretization is carried out by choosing 
$N\in \mathbb{N}$, the number of time steps.
To simplify the discussion, we assume uniform time steps, let $\Delta t = \mathbb{T}/N$.
We set $t^k = k \Delta t$ and for a time dependent function we denote $\varphi^k = \varphi(t^k)$.
Over these sequences we define the operators
\begin{equation}
\label{eq:defofbdf}
\bdf{\varphi^{k+1}}{m} := 
\begin{cases}
\frac{1}{\Delta t}(\varphi^{k+1}-\varphi^{k}) & m = 1, \\
\frac{1}{2\Delta t} \left( 3 \varphi^{k+1} - 4\varphi^k + \varphi^{k-1} \right) & m = 2,
\end{cases}
\end{equation}
for the backward Euler time discretization order 1 and order 2.
In this paper, we employ BDF2 (second order backward difference formula) with $m=2$ to discretize the time derivatives.

Thus we obtain the following time discretized formulation 
\begin{multline}
\label{eqn:slight_comp_simp_time_1}
\phi c_w^F \rho_w^0  \bdf{s_w^{k+1} p^{k+1}_w}{m} 
-\nabla \cdot \left(\bK (\lambda_t(s_w^{k+1}) \nabla p_w^{k+1} 
 - (\rho^0 \lambda(s_w^{k+1}))_t \bg ) \right) \\
 -\nabla \cdot \left( \bK \lambda_n(s_w^{k+1})  \nabla p_c(s_w^{k+1}) \right)  = (\rho^0 f^{k+1})_t,
\end{multline}
As frequently done in modeling slightly compressible two-phase flow, we neglect the terms involving small compressibility $c_n^F$ in \eqref{eqn:slight_comp} with the exception of  $c_w^F$. 
Here  $c_w^F$ is included as a regularization term for the solver.  

Next, the saturation system is discretized by 
\begin{equation}
\phi \rho^0_w \bdf{s^{k+1}_w}{m}
 + \nabla \cdot 
\left(-\rho^0_w \bK \dfrac{k_w(s_w^{k+1}) }{\mu_w} ( \nabla p_w^{k+1} - \rho_w \bg ) \right) 
 = \rho^0_w f_w^{k+1},
 \label{eqn:transport_time_1}
\end{equation}
The above system is fully coupled and nonlinear. We propose  the  following iterative decoupled scheme.

\subsubsection{Sequential IMPES algorithm}
\label{sec:IMPES}

The implicit pressure and explicit saturation algorithm (IMPES) is frequently applied as an efficient algorithm for decoupling 
and sequentially solving the system \cite{ChenHuan2006}. 
For uniform time steps, to approximate the time dependent terms we define the extrapolation of 
$\varphi^{k+1,*}$ by
$$
\varphi^{k+1,*} :=
\varphi^k + ( \varphi^k - \varphi^{k-1}).
$$
The IMPES algorithm solves the system as follows:
\begin{enumerate}
\item Initial conditions  at time $t^{k-1}, t^k$ are given. 
\item  Solve $p^{k+1}_w$ at time $t^{k+1}$ by using the previous  saturation to compute $\lambda_i(s_w^{k+1,*})$ and $p_c(s_w^{k+1,*})$. 
\begin{multline}
\label{eqn:incomp_time_impes}
\phi c_w^F \rho_w^0  \bdf{s_w^{k+1,*} p^{k+1}_w}{m} 
-\nabla \cdot \left(\bK \lambda_t(s_w^{k+1,*}) \nabla p_w^{k+1}\right)   \\
 = (\rho^0 f)_t
 -\nabla \cdot \left(\bK 
 (\rho^0 \lambda(s_w^{k+1,*}))_t \bg ) \right)  
+\nabla \cdot \left(  \bK \lambda_n(s_w^{k+1,*}) \nabla p_c(s_w^{k+1,*}) 
\right) 
\end{multline}
\item Compute the velocity $\bu_w^{k+1,*}$ by using $p_w^{k+1}$ and the saturation. 
\item Compute $s_w^{k+1}$ using an explicit time stepping.
\begin{equation}
\phi \rho^0_w \bdf{s^{k+1}_w}{m}
 = \rho^0_w f_w^{k+1}
  + \nabla \cdot 
 \left(\rho^0_w \bK   \dfrac{k_w(s_w^{k+1,*}) }{\mu_w} \left( \nabla p_w^{k+1} - \rho^0_w \bg \right) \right) 
 \label{eqn:transport_time_impes}
\end{equation}

\end{enumerate}

\subsubsection{Iterative IMPES algorithm}
\label{sec:itIMPES}
An iterative IMPES algorithm is to  solve the following equations sequentially for iterations $j=1,\cdots$ until it converges to a given tolerance or a fixed number of iterations has been reached.
For example, at each time step $t^k$: 

\begin{enumerate}
\item For $j=0$, set $s_w^{k+1,j} = s_w^k$ and  $s_w^{k+1,j-1} = s_w^{k-1}$. 
Solve for $p_w^{k+1,j+1}$  satisfying 
\begin{multline}
\label{eqn:incomp_time_impes_it}
\phi c_w^F \rho_w^0 \bdf{s_w^{k+1,*,j} p^{k+1}_w}{m} 
-\nabla \cdot \left(\bK \lambda_t(s_w^{k+1,*,j}) \nabla p_w^{k+1,j+1}\right) \\
 = (\rho^0 f)_t
 -\nabla \cdot \left(\bK 
 (\rho^0 \lambda(s_w^{k+1,*,j}))_t \bg ) \right)  
+\nabla \cdot \left(  \bK \lambda_n(s_w^{k+1,*,j}) \nabla p_c(s_w^{k+1,*,j}) 
\right) ,
\end{multline}
where $s_w^{k+1,*,j}  = s_w^{k+1,j}  + (s_w^{k+1,j} - s_w^{k+1,j-1})$.
\item Given $s_w^{k+1,*,j}$ and $p_w^{k+1,j+1}$, solve for $s_w^{k+1,j+1}$ satisfying
\begin{equation}
\phi \rho^0_w \bdf{s^{k+1,j+1}_w}{m}
 = \rho^0_w f_w^{k+1}
  + \nabla \cdot 
 \left(\rho^0_w \bK  \dfrac{k_w(s_w^{k+1,*,j}) }{\mu_w}  \left( \nabla p_w^{k+1,j+1} - \rho^0_w \bg \right) \right). 
 \label{eqn:transport_time_impes_it}
\end{equation}
\item Iteration continues until $\| s_w^{k+1,j+1} - s_w^{k+1,j} \| \leq \varepsilon_{I}$. 
\end{enumerate}

\subsection{Spatial Approximation of the Pressure System}
The locally conservative EG is selected
for the space approximation of the pressure system  \eqref{eqn:slight_comp_simp_time_1}.
Here we apply the discontinuous Galerkin (DG)
IIPG (incomplete interior penalty Galerkin) method 
for the flow problem to satisfy the discrete sum compatibility condition 
\cite{DawsonC_SunS_WheelerM-2004aa, LeeWhe_2017, ScovazziWheMikLee16}.
Mathematical stability and error convergence of EG for a single phase system
is discussed in \cite{LeeLeeWhi15,LeeLeeWhi17a,LeeWhe_2017}

The EG finite element space approximation of the wetting phase pressure $p_w(\bx,t)$ is denoted by $P_w(\bx,t) \in V^{\textsf{EG}}_{h,l}(\mathcal{T}_h) $ and 
we let $P_w^k := P_w(\bx,t^k)$ for time discretization, 
$0 \leq k \leq N$.
We set an initial condition for the pressure as  $P_w^0 := \Pi_h  p_w(\cdot,0)$.
Let  $p_\textsf{in}^{k+1}, p_\textsf{out}^{k+1}, \bu_N^{k+1}$ and $f^{k+1}$ are approximations of $p_\textsf{in}(\cdot, t^{k+1}), p_\textsf{out}(\cdot, t^{k+1}), \bu_N(\cdot, t^{k+1})$ and $f(\cdot, t^{k+1})$ on $\Gamma_D$, $\Gamma_N$ and $\Omega$, respectively at time $t^{k+1}$.
Assuming  $s_w(\cdot,t^{k+1})$  is known, 
and employing time lagged/extrapolated values for simplicity, 
the time stepping algorithm reads as follows: Given $P_w^{k-1}$, $P_w^{k}$, find
\begin{equation}\label{eq:egscheme}
P_w^{k+1} \in V_{h,l}^{\textsf{EG}}(\mathcal{T}_h) \mbox{  such that  } \calS(P_w^{k+1},\omega) = \mathcal{F}(\omega), \quad \forall\, \omega \in V_{h,l}^{\textsf{EG}}(\mathcal{T}_h) , \,
\end{equation}
where $\calS$ and $\mathcal{F}$ are the bilinear form and linear functional, respectively, are defined as 
\begin{multline*}
\calS(P_w^{k+1},\omega) := 
\left( (\phi \rho_w^0  c_w^F s_w^{k+1,*} 
)\dfrac{ 3 }{2 \Delta t}P_w^{k+1}, \omega \right)_{\mathcal{T}_h} 
+ 
\left (\lambda_t(s_w^{k+1,*}) \bK \nabla P_w^{k+1},\nabla \omega \right )_{\Th} \\
- \left \langle  \av{\bK \lambda_t(s_w^{k+1,*})  \nabla P_w^{k+1} }, \jump{\omega} \right \rangle_{\mathcal{E}_h^{1}}  
+ {\dfrac{\alpha}{h_{e}} }\av{\bK  \lambda_t(s_w^{k+1,*})} \left\langle    \jump{P_w^{k+1}},\jump{\omega} \right\rangle_{\mathcal{E}_h^{1}},
\end{multline*}
and 
\begin{multline*}
\mathcal{F}(\omega) := 
\left( 
(\phi \rho^0_w  c^F_w s_w^{k} )
(\dfrac{ 2 }{\Delta t}P_w^{k})
-
(\phi \rho^0_w  c^F_w s_w^{k-1} )
(\dfrac{1}{2 \Delta t}P_w^{k-1} )
, \omega 
\right)_{\mathcal{T}_h} 
+
\left( (\rho^0 f^{k+1})_t,\omega\right)_{\mathcal{T}_h}  \\
- \left(\bK \lambda_n(s_w^{k+1,*})  \nabla p_c(s_w^{k+1,*}) 
-  \bK(\rho^0 \lambda(s_w^{k+1,*}))_t     \bg, \nabla \omega
\right)_{\mathcal{T}_h} 
+
\left<
\av{\bK \lambda_n(s_w^{k+1,*})  \nabla p_c(s_w^{k+1,*}) 
-\bK  (\rho^0  \lambda (s_w^{k+1,*}))_t   \bg  },\jump{\omega}
\right>_{\mathcal{E}_h^{1}} \\
{- \dfrac{\alpha_c}{h_{e}}\av{ \bK  \lambda_n(s_w^{k+1,*})} \left\langle    \jump{p_c(s_w^{k+1,*})},\jump{\omega} \right\rangle_{\mathcal{E}_h^{1}},} \\
+ 
{\dfrac{\alpha}{h_{e}}}\av{ \bK  \lambda_t(s_w^{k+1,*})}  \left\langle  p_\textsf{in}^{k+1},\jump{\omega}\right\rangle_{\mathcal{E}_h^{{\textsf{in}},\partial}}
+ 
{\dfrac{\alpha}{h_{e}}}\av{ \bK  \lambda_t(s_w^{k+1,*})}  \left\langle  p_\textsf{out}^{k+1},\jump{\omega}\right\rangle_{\mathcal{E}_h^{{\textsf{out}},\partial}}
- \left \langle  {{\bu}^{k+1}_N},\jump{\omega} \right \rangle_{\mathcal{E}_h^{N,\partial}}.
\end{multline*}

Here 
$h_e$ denotes the maximum {length} of the edge $e \in \Eh$ and 
$\alpha, \alpha_c$ are penalty parameters for pressure and capillary pressure, respectively. 
For adaptive mesh refinement with hanging nodes, we make the usual assumption to 
set the $h_{e} = \min(h^+,h^-)$ for $e=\partial \K^{+}\cap \partial \K^{-}$ over the edges on a mesh $T$.

\subsubsection{Locally conservative flux}
\label{sec:flux}
Conservative flux variables are described in \cite{LeeLeeWhi15,sunliu2009} 
with details for convergence analyses.
With slight modifications to the latter single phase case, we define the 
two-phase wetting phase velocity as
$\bU_w^{k+1,*}$ since it depends on the previous saturation value $s_w^{k+1,*}$.
Let $P_w^{k+1}$ be the wetting phase solution to \eqref{eq:egscheme}, then  we define the globally and locally conservative flux variables $\bU_w^{k+1,*}$ at time step $t^{k+1}$ by the following :
\begin{align}\label{flux}
\bU_w^{k+1,*} |_{T} &:= -\bK\dfrac{k_w(s_w^{k+1,*})}{\mu_w}  \left( \nabla P_w^{k+1} - \rho_w^0 \bg \right), \;  \forall T \in \mathcal{T}_h  \vspace*{0.2in} \\ 
\bU_w^{k+1,*} \cdot \bn|_{e} &:= - 
 \av{ \bK\dfrac{k_w(s_w^{k+1,*})}{\mu_w}  \left( \nabla  P_w^{k+1} - \rho_w^0 \bg \right)} \cdot \bn 
+ \av{\dfrac{\alpha}{h_e} \bK \dfrac{k_w(s_w^{k+1,*})}{\mu_w}  }
 \jump{P_w^{k+1}}, 
\; \forall e \in \mathcal{E}_h^I, \\ 
\bU_w^{k+1,*} \cdot \bn|_{e} &:= \bu_{Nw}^{k+1}, \; \forall e \in \mathcal{E}_h^{N,\partial},  \\ 
\bU_w^{k+1,*} \cdot \bn|_{e} &:= - \bK \dfrac{k_w(s_w^{k+1,*})}{\mu_w} \left(\nabla P_w^{k+1} - \rho_w^0 \bg \right) \cdot \bn 
+  \av{ \dfrac{\alpha}{h_e}\bK  \dfrac{ k_w(s_w^{k+1,*})}{\mu_w}  }
  \left( P_w^{k+1} - p_\textsf{in/out}^{k+1} \right ), \; \forall e \in \mathcal{E}_{h_e}^{D,\partial},
\label{flux-end}
\end{align}
where $\bn$ is the unit normal vector of the boundary edge $e$ of $T$ and
$\bu_{Nw}^{k+1} := (\bu_{N}^{k+1} - \bK \lambda^{k+1}_n\nabla p_c^{k+1})(\lambda_w^{k+1}/(\lambda_t^{k+1} \rho^0_w))$.

\subsubsection{H(div) reconstruction of the flux}
\label{sec:hdiv}
For incompressible flow, it is frequently useful to project  the velocity (flux) into a 
$H$(div) space for high order approximation to a transport system, see  \cite{Arbogast:2013eg, bastian2003superconvergence,ern2009accurate, 
Ern:2010cn,ern2007accurate} for more details.
We illustrate below, the reconstruction of the EG flux \eqref{flux}-\eqref{flux-end} in a $H$(div) space for quadrilateral elements \cite{LeeLeeWhi17a,Li:2015fl}.
The flux is projected into the Raviart-Thomas  (RT$_l$) space \cite{boffi2013mixed,Raviart1977},
$$
\mathcal{H} := \{ v \in H(\text{div}) : 
v|_E \in \mathbb{Q}_{l+l,l}(T) \times \mathbb{Q}_{l,l+1}(T) , 
\ \forall T \in \mathcal{T}_h \}, 
$$
where 
$$
\mathbb{Q}_{a,b}(T) := 
\left\{ v \ : \ v(\bx) = \sum_{i=0}^a\sum_{j=0}^b \omega_{i,j} \bx_1^i\bx_2^j, \ \bx \in T, \omega_{i,j} \in \mathbb{R} \right\}
$$
with polynomial order $l$.

Let ${\bU}^{div} \in \mathcal{H}$ be the reconstructed flux defined on each element $T$ as 
\begin{equation}
({\bU}^{div} , v)_T 
= (\bU, v)_T,
\end{equation} 
where $v \in \mathbb{Q}_{l-1,l}(T) \times \mathbb{Q}_{l,l-1}(T)$ 
and
\begin{equation}
\langle {\bU}^{div}\cdot \bn , w \rangle_e 
= \langle  \bU \cdot \bn, w \rangle_e, \ \forall e \in \partial T, \
w \in \mathbb{P}_l(e).
\end{equation}
We note that the polynomial order of the post-processed space $\mathcal{H}$ is chosen consistently with the order of the pressure space $l$.
The performance of the projection is illustrated in \cite{LeeLeeWhi17a}.

\subsection{Spatial Approximation of the Saturation System}
The bilinear form of EG coupled with an entropy residual stabilization is employed for modeling the transport system \eqref{eqn:transport} with  high order approximations \cite{LeeWhe_2017}. 
Here, again we apply DG IIPG method although other interior penalty methods can be utilized. 
Stability and error convergence analyses for the approximation are provided in \cite{LeeLeeWhi17a}.

The EG finite element space approximation of the wetting phase saturation $s_w(\bx,t)$ is denoted by $S_w(\bx,t) \in V^{\textsf{EG}}_{h,s}(\mathcal{T}_h) $ and 
we let $S_w^k := S_w(\bx,t^k)$ for time discretization, $0 \leq k \leq N$.
We set an initial condition for the saturation as  
$S_w^0 := \Pi_h  s_w(\cdot,0)$.
With $P_w^{k+1}$ computed by the system \eqref{eq:egscheme} and locally conservative fluxes \eqref{flux}, 
the time stepping algorithm reads as follows: Given $S_w^{k-1}$,$S_w^{k}$, find
\begin{equation}\label{eq:egscheme_saturation}
S_w^{k+1} \in V_{h,{s}}^{\textsf{EG}}(\mathcal{T}_h) \mbox{  such that  } \calM(S_w^{k+1},\psi) = \mathcal{G}(\psi), \quad \forall\, \psi \in V_{h,s}^{\textsf{EG}}(\mathcal{T}_h) , \,
\end{equation}
where, 
\begin{equation}
\calM(S_w^{k+1},\psi) = 
\left(  \phi \rho_w^{0}\dfrac{3}{2\Delta t} S_w^{k+1}, \psi  \right)_{\mathcal{T}_h}
{- (\rho_w^{0} S_w^{k+1}  (f_w^{k+1})^-,  \psi)_{\mathcal{T}_h} }
\end{equation}
and 
\begin{align}
\mathcal{G}(\psi) 
&= 
\left(  \phi \dfrac{2\rho_w^{0}}{\Delta t} S_w^{k}
-\phi \dfrac{\rho_w^{0}}{2\Delta t} S_w^{k-1}
, \psi \right)_{\mathcal{T}_h}
{+( \rho_w^{0} (f^{k+1}_w)^+, \psi)_{\mathcal{T}_h}}
- ({\rho_w^{0}} \nabla \cdot \bU_w^{k+1,*}, \psi)_{\mathcal{T}_h} \nonumber \\
&=
\left(  \phi \dfrac{2\rho_w^{0}}{\Delta t} S_w^{k}
-\phi \dfrac{\rho_w^{0}}{2\Delta t} S_w^{k-1}, \psi \right)_{\mathcal{T}_h} 
+(\rho_w^{0} { (f^{k+1}_w)^+ }, \psi)_{\mathcal{T}_h}
+ ({\rho_w^{0}}\bU_w^{k+1,*}, \nabla \psi)_{\mathcal{T}_h} 
- \left< {\rho_w^{0}} \bU_w^{k+1,*} \cdot \bn, \jump{\psi} \right>_{\Eh}
\end{align}
The injection/production term $f_w^{k+1}:= \tilde{s}_w^{k+1} q^{k+1}_w $ splits by 
$$
(f_w^{k+1})^+ = \max(0,f_w^{k+1}) \quad \text{ and } \quad 
(f_w^{k+1})^- = \min(0,f_w^{k+1}).
$$
Recall that $\tilde{s}_w^{k+1}$ is the injected saturation if $q_w^{k+1} > 0$ 
and is the resident saturation if $q_w^{k+1} < 0$. 
The computed locally conservative numerical fluxes in the section \ref{sec:flux} are applied here.

\subsubsection{Entropy residual stabilization}
\label{sec:levelset_entropy}
Elimination of spurious numerical oscillations due to sharp gradients in the solution requires stabilizations for the high order approximation to the transport system ($s \geq 1$).
In this section, we describe an entropy viscosity  stabilization technique 
to avoid oscillations in the EG formulation \eqref{eq:egscheme_saturation}.
This method was introduced in \cite{guermond2011entropy} and mathematical stability properties are discussed in 
\cite{MR3167449} for CG and in \cite{Zingan:2013bb} for DG. Recently, it was employed for EG single phase miscible displacement problems \cite{LeeWhe_2017} by the authors. 
Here, we provide an extension to two-phase flow saturation equation. 

We redefine the velocity term for the two-phase flow system by  separating the relative permeability which is a function of saturation, as is frequently 
referred to as expanded mixed form \cite{arbogast1997mixed}. We let
\begin{align}
\bu_i &= - \bK \dfrac{k_i(s_w) }{\mu_i} \left( \nabla p_i - \rho_i \bg \right) \\
&= k_i(s_w)  \hat{\bu}_i,
\end{align}
where 
\begin{equation}
\hat{\bu}_i := - \dfrac{\bK}{\mu_i} \left( \nabla p_i - \rho_i \bg \right), \; \; i\in \{n,w\}.
\label{eqn:new_vel}
\end{equation}

Now, we introduce a numerical dissipation term 
${\mathcal{E}}(S_w^{k+1}, \psi)$ in  \eqref{eq:egscheme_saturation} to obtain, 
\begin{equation}
\mathcal{M}(S_w^{k+1}, \psi) + 
\mathcal{E}(S_w^{k+1}, \psi)
= {\mathcal{G}}(\psi), \quad \forall \psi  \in V_{h,s}^{\textsf{EG}}(\mathcal{T}_h),
\label{eqn:discrete_transport_ev} 
\end{equation}
where 
\begin{multline}
{\mathcal{E}}(S_w^{k+1}, \psi)  
:= 
\left( {\rho^0_w}  {\mu}^{k+1}_{\text{Stab}}(S_w,\hat{\bU}_{i}) _{|T} \nabla S_w^{k+1}, \nabla \psi \right)_{\mathcal{T}_h}  \\
-
\left\langle \av {\rho^0_w {\mu}^{k+1}_{\text{Stab}}(S_w,\hat{\bU}_{i}) _{|T}  \nabla S_w^{k+1} }, \jump{\psi} \right\rangle_{\mathcal{E}_h^{I} }
+
 \left\langle \av{ \dfrac{\alpha_T}{h_{e}} \rho^0_w {\mu}^{k+1}_{\text{Stab}}(S_w,\hat{\bU}_{i}) _{|T}}  \jump{S_w^{k+1}}, \jump{\psi} \right\rangle_{\mathcal{E}_h^{I}},
\end{multline}
and $\alpha_T$ is a penalty parameter.

Here ${\mu}^{k+1}_{\text{Stab}}(S_w,\hat{\bU}_{i}) _{|T}: 
\Omega \times [0,\mathbb{T}] \rightarrow \mathbb{R}$ is
the
stabilization coefficient, which is piecewise constant over the mesh $T$. 
It is defined on each $T \in \mathcal{T}_h$ by 
\begin{equation}\label{eqn:levelset_stab_v}
\rho^0_w{\mu}^{k+1}_{\text{Stab}}(S_w,\hat{\bU}_{i}) _{|T} := 
\min(
\rho^0_w \mu^{k+1}_{{\textsf{Lin}}}(S_w,\hat{\bU}_{i})_{|T} , 
\rho^0_w \mu^{k+1}_{{\textsf{Ent}}}(S_w,\hat{\bU}_w)_{|T} ).
\end{equation}
The main idea of the entropy residual stabilization is to split the 
stabilization terms into $\mu^{k+1}_{{\textsf{Lin}}}$ and $\mu^{k+1}_{{\textsf{Ent}}}$.
If $S_w(\cdot, t)$ is smooth, the entropy viscosity stabilization
$\mu^{k+1}_{{\textsf{Ent}}}(S_w,\hat{\bU}_w)_{|T}$ will be activated, 
since $\mu^{k+1}_{{\textsf{Ent}}}$ is small.
However, 
the linear viscosity ${\mu}^{k+1}_{{\textsf{Lin}}}(S_w,\hat{\bU}_{i})_{|T}$ is activated 
where $S_w(\cdot, t)$ is not smooth.
The first order linear viscosity is defined by, 
\begin{equation}
{\mu}^{k+1}_{{\textsf{Lin}}}(S_w,\hat{\bU}_w)_{|T}  :=  \lambda_{{\textsf{Lin}}} h_T 
\|\max_{i\in\{n,w\}}(k'_i({S^{k+1,*}_w}) \ \hat{\bU}_{i}^{k+1}) \|_{L^{\infty}(T)}, 
\quad \forall T \in \mathcal{T}_h ,
\label{eqn:level:fv}
\end{equation}
where  $h_T$ is the mesh size and 
$\lambda_{{\textsf{Lin}}}$ is a positive constant.
We note that $s_w$ is transported by $\hat{\bu}_w$ and 
$s_n = 1-s_w$ is transported by $\hat{\bu}_n$.

Next, we describe the entropy viscosity stabilization. 
Recall that it is known that the scalar-valued conservation equation
\begin{equation}
{\partial_t} (\phi \rho_w s_w) + \nabla \cdot  
{\boldsymbol{v}}(s_w) = \rho_w  {f_w}
\end{equation}
may have one weak solution 
in the sense of distributions satisfying the additional inequality 
\begin{equation}
{\partial_t} (\phi \rho_w E(s_w)) + 
\nabla \cdot 
{\boldsymbol{F}}(s_w)
- E'(s_w) \rho_w {f_w} \leq 0, 
\end{equation}
for any convex function $E \in \mathcal{C}^0(\Omega;\mathbb{R})$ which is called entropy and 
${\boldsymbol{F}}'(s_w):= E'(s_w) {\boldsymbol{v}}'(s_w)$, the associated entropy flux \cite{Kruzkov,Panov1994}. 
The equality holds for smooth solutions. 

For the two-phase flow system, we redefined the velocity in \eqref{eqn:new_vel} to split the relative permeability. Thus, we set
${\boldsymbol{v}}(s_w) :=  {\rho^0_w} k_w(s_w) \hat{\bu}_w$.
Then we obtain 
${\boldsymbol{F}}'(s_w) = ({\rho^0_w} k_w'(s_w) \hat{\bu}_w) \cdot E'(s_w)$ and 
$\nabla \cdot {\boldsymbol{F}}(s_w) = {\boldsymbol{F}}'(s_w) \cdot \nabla s_w$.
Note that we can rewrite $\nabla E(s_w) = E'(s_w)\nabla s_w$.
We define the entropy residual which is a reliable indicator of the regularity of $s_w$ as 
\begin{equation}
R_{\textsf{Ent}}^{k+1}(S_w,\hat{\bU}_w) := 
\bdf{\phi \rho_w E(S_w^k)}{m}  
+ {\rho^0_w} k_w'({S_w^{k+1,*}}){\hat{\bU}_w^{k+1}} E'(S_w^{k+1,*})\nabla (S_w^{k+1,*}) - E'(S_w^{k+1,*}) \rho_w f_w,
\label{eqn:ent_res_1}
\end{equation}
which is large when $S_w$ is not smooth. 
In this paper, we chose 
\begin{equation}
E(S_w^{k+1,*}) = \dfrac{1}{b} |S_w^{k+1,*}|^b, \ b  \text{ is a positive even number}  
\label{eqn:entropy_func_1}
\end{equation}
with $b=10$ or 
\begin{equation}
E(S_w^{k+1,*}) = - \log( | S_w^{k+1,*} (1-S_w^{k+1,*})| + \varepsilon)
\label{eqn:entropy_func_2}
\end{equation}
with $\varepsilon < 1$ as chosen in \cite{FLD:FLD4071,Guermond2017448,LeeWhe_2017}.
Finally, the local entropy viscosity for each step is defined as
\begin{equation}
\mu^{k+1}_{{\textsf{Ent}}}(S_w, \hat{\bU}_w)_{|T} := \lambda_{\textsf{Ent}} h_T^2 
\dfrac{{ER_{\textsf{Ent}}^{k+1}}_{|T} } 
{\|E(S_w^{k+1,*}) - \bar{E}^{k+1,*} \|_{L^{\infty}(\Omega)}}, 
\quad  \forall T \in \mathcal{T}_h,
\label{eqn:ent_vis}
\end{equation}
where 
\begin{equation}
{ER_{\textsf{Ent}}^{k+1}}_{|T} := 
\max(\|R_{\textsf{Ent}}^{k+1}\|_{L^{\infty}(T)}, \|J_{\textsf{Ent}}^{k+1}\|_{L^{\infty}(\partial T)}).
\label{eqn:ent_res}
\end{equation}
Here $\lambda_{\textsf{Ent}}$ is a positive constant to be chosen with the average $\bar{E}^{k+1,*} := \frac{1}{|\Omega|} \int_{\Omega}E(S_w^{k+1,*}) \ d\bx$. 
We define the residual term calculated on the faces by
\begin{equation}
J_{\textsf{Ent}}^{k+1}(S_w,\hat{\bU}_w) := 
h^{-1}_{T}  \av{\hat{\bU}_w^{k+1}} \cdot \jump{E(S_w^{k+1,*})}.
\end{equation}
The entropy stability with above residuals for discontinuous case is given 
with  more details
in \cite{Zingan:2013bb}. 
Also, readers are referred to \cite{guermond2011entropy} for tuning the constants 
($\lambda_{\textsf{Ent}}, \lambda_{\textsf{Lin}}$). 

%
%
\subsection{Adaptive Mesh Refinement}
In this section, we propose a refinement strategy by increasing the mesh resolution in the cells where the entropy residual values \eqref{eqn:ent_res} are locally larger than others. It is shown in \cite{andrews1998posteriori,puppo2004numerical} that the entropy residual can be used as a posteriori error indicator.
The general residual of the system \eqref{eq:egscheme_saturation} 
could also be utilized as an error indicator, but this residual goes to zero as $h \rightarrow 0$ due to consistency.
However, as discussed in \cite{guermond2011entropy}, 
the entropy residual \eqref{eqn:ent_res} converges to a Dirac measure supported in the neighborhood of shocks. 
In this sense, the entropy residual is a robust indicator 
and also efficient 
since it is  been computed for a stabilization.

\begin{figure}[!h]
\centering
\begin{tikzpicture}[scale=1.25]
\node[] at (1,1) {\footnotesize $\textsf{Ref}_T=0$};
\node[] at (4.5,1.5) {\footnotesize$\textsf{Ref}_T=1$};
\node[] at (7.5,1.75) {\footnotesize$\textsf{Ref}_T=2$};

\draw (0,0) -- (2,0); 
\draw (2,0) -- (2,2); 
\draw (2,2) -- (0,2); 
\draw (0,2) -- (0,0); 

\node[] at (2.5,1) {$\Rightarrow$} ; 

\draw (3,0) -- (5,0); 
\draw (5,0) -- (5,2); 
\draw (5,2) -- (3,2); 
\draw (3,2) -- (3,0); 
\draw (4,0) -- (4,2); 
\draw (3,1) -- (5,1); 

\node[] at (5.5,1) {$\Rightarrow$} ; 

\draw (6,0) -- (8,0); 
\draw (8,0) -- (8,2); 
\draw (8,2) -- (6,2); 
\draw (6,2) -- (6,0); 
\draw (7,0) -- (7,2); 
\draw (6,1) -- (8,1); 
\draw (7.5,1) -- (7.5,2); 
\draw (7,1.5) -- (8,1.5); 

\draw[solid] (7.5,1) circle (0.05);
\draw[solid] (7.,1.5) circle (0.05);
\end{tikzpicture}
\caption{Adaptive mesh refinement levels. $
\textsf{Ref}_T$ is the refinement level and $\circ$ denotes the hanging nodes. The mesh refines until $\textsf{Ref}_T < R_{\max}$.}
\label{fig:mesh}
\end{figure}
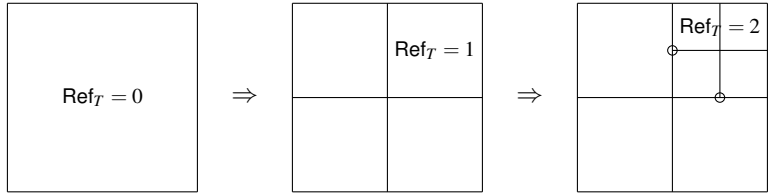

We denote the refinement level, $\textsf{Ref}_T$ (see Figure \ref{fig:mesh}), to be  the number of times a cell($T$) from the initial subdivision has been refined to produce the current cell. 
Here, a cell $T$ is refined if its corresponding $\textsf{Ref}_T$ is smaller than a given number $R_{\max}$
 and if 
\begin{equation}
|{ER_{\textsf{Ent}}^{k+1}}_{|T}(\bx_T,t)|  \geq 
C_R 
\max_{\K \in \Th} |{ER_{\textsf{Ent}}^{k+1}}_{|T}(\bx_T,t)|,
\end{equation}
where $\bx_T$ is the barycenter of $T$ and $C_R \in [0,1]$. The purpose of the parameter $R_{\max}$ is to control the total number of cells, which is set to be two more than the initial 
$\textsf{Ref}_T$.
A cell $T$ is coarsened if 
\begin{equation}
|{ER_{\textsf{Ent}}^{k+1}}_{|T}(\bx_T,t)| \leq C_C
\max_{\K \in \Th} |{ER_{\textsf{Ent}}^{k+1}}_{|T}(\bx_T,t)|,
\end{equation}
where
$C_C \in [0,1]$.
However, a cell is not coarsened if the 
$\textsf{Ref}_T$ is smaller than a given number $R_{\min}$.
Here $R_{\min}$ is set to be two less than the initial $\textsf{Ref}_T$.
In addition, a cell is not refined more if the total number of cells are more than $\textsf{Cell}_{\max}$.
The subdivisions are accomplished with at most one hanging node per face.
During mesh refinement, to initialize or remove  nodal values,
standard interpolations and restrictions are employed, respectively. 
We take advantage of the dynamic mesh adaptivity feature with hanging nodes in deal.II \cite{dealII84} in which subdivision and mesh distribution are implemented using the p4est library \cite{p4est}.

\subsection{Global Algorithm and Solvers}
\tikzstyle{block00} = [rectangle, draw, text centered, 
rounded corners, text width= 4em,   minimum height=4.em, node distance= 2.cm]
\tikzstyle{block10} = [rectangle, draw, text centered, 
rounded corners, text width= 4em,   minimum height=4.em, node distance= 2.2cm]
\tikzstyle{block11} = [rectangle, draw, text centered, 
rounded corners, text width= 6em,   minimum height=4em, node distance= 2.5cm]
\tikzstyle{block12} = [rectangle, draw, text centered, 
rounded corners, text width=6em,    minimum height=4em, node distance= 2.8cm]
\tikzstyle{block13} = [rectangle, draw, text centered, 
rounded corners, text width=4em,    minimum height=4em, node distance= 2.45cm]
\tikzstyle{line11}  = [draw, -latex']
\tikzstyle{line} = [draw, -latex']
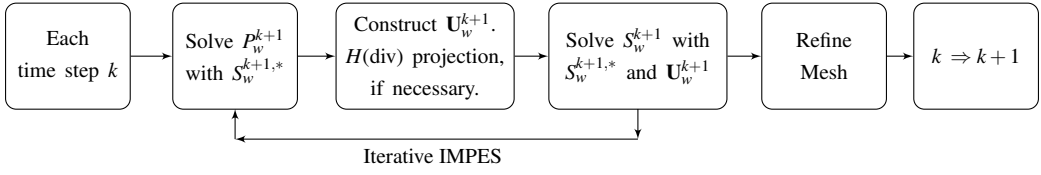
\begin{figure}[!h]
\centering
\begin{tikzpicture}[scale=1., every node/.style={transform shape}]
\node [block00] (Init) {\footnotesize Each \\ time step $k$};
\node [block10, right of=Init] (sol_vel) {\footnotesize 
               Solve $P_w^{k+1}$ \\ 
               with $S_w^{k+1,*}$};
\node [block11, right of=sol_vel] (sol_pres) {\footnotesize 
        Construct $\bU_w^{k+1}$. $H$(div) projection, if necessary.};
\node [block12, right of=sol_pres] (sol_sat) {\footnotesize 
        Solve $S_w^{k+1}$ with $S_w^{k+1,*}$ and $\bU_w^{k+1}$};
\node [block13, right of=sol_sat] (end) { \footnotesize Refine Mesh};
\node [block00, right of=end] (Init1) {\footnotesize $k$ $\Rightarrow k+1$};

\path [line11] (Init) -- (sol_vel);
\path [line11] (sol_vel) -- (sol_pres);
\path [line11] (sol_pres) -- (sol_sat);
\path [line11] (sol_sat) -- (end);
\path [line11] (end) -- (Init1);

\draw[] (4.8,-1.1) node[below]  {\footnotesize Iterative IMPES};
\path [line] (7.5,-0.7) -- (7.5,-1.1); 
\path [line] (7.5,-1.1) -- (2.2, -1.1); 
\path [line] (2.2, -1.1) -- (2.2,-0.7); 

\end{tikzpicture}
\caption{Flowchart of global solution algorithm.}
\label{fig:arg_flow}
\end{figure}

We present our global algorithm in Figure \ref{fig:arg_flow} for modeling the two-phase flow  problem.  
An efficient solver developed in \cite{LeeLeeWhi15} is applied to
solve the EG pressure and saturation system separately. The current solver is GMRES  Algebraic Multigrid(AMG) block diagonal preconditioner.
$H$(div) projection is activated only for incompressible cases.
The entropy residuals are employed when solving the transport system as well as refining the mesh.
The authors created the EG two-phase flow code to compute the following numerical examples based on the open-source finite element package deal.II \cite{dealII84} which is coupled with the parallel MPI library \cite{open_mpi} and Trilinos solver \cite{Trilinos-Overview}.

\section{Numerical Examples}
\label{sec:num_ex}
This section verifies and demonstrates the performance of our proposed EG algorithm. 
First, 
the convergence of the spatial errors  are shown for the two-phase EG flow system for decoupled, sequential and iterative IMPES.
Next, several numerical examples with 
capillary pressure, gravity and dynamic mesh adaptivity including a benchmark test are provided.

\subsection{Example 1. Convergence Tests - decoupled case with entropy residual stabilization.}

Here we consider the two-phase flow problem with exact solution given by 
\begin{equation}
p_w = \cos(t+x-y), \; \; \; s_w = \sin(t+x-y+1)
\end{equation}
in the domain $\Omega = (0,\SI{1}{\metre})^2$.
A Dirichlet boundary condition is applied for the pressure system.
\begin{figure}[!h]
\centering
\subfloat[Capillary pressure with epsilon]
{
\includegraphics[width=0.3\textwidth]{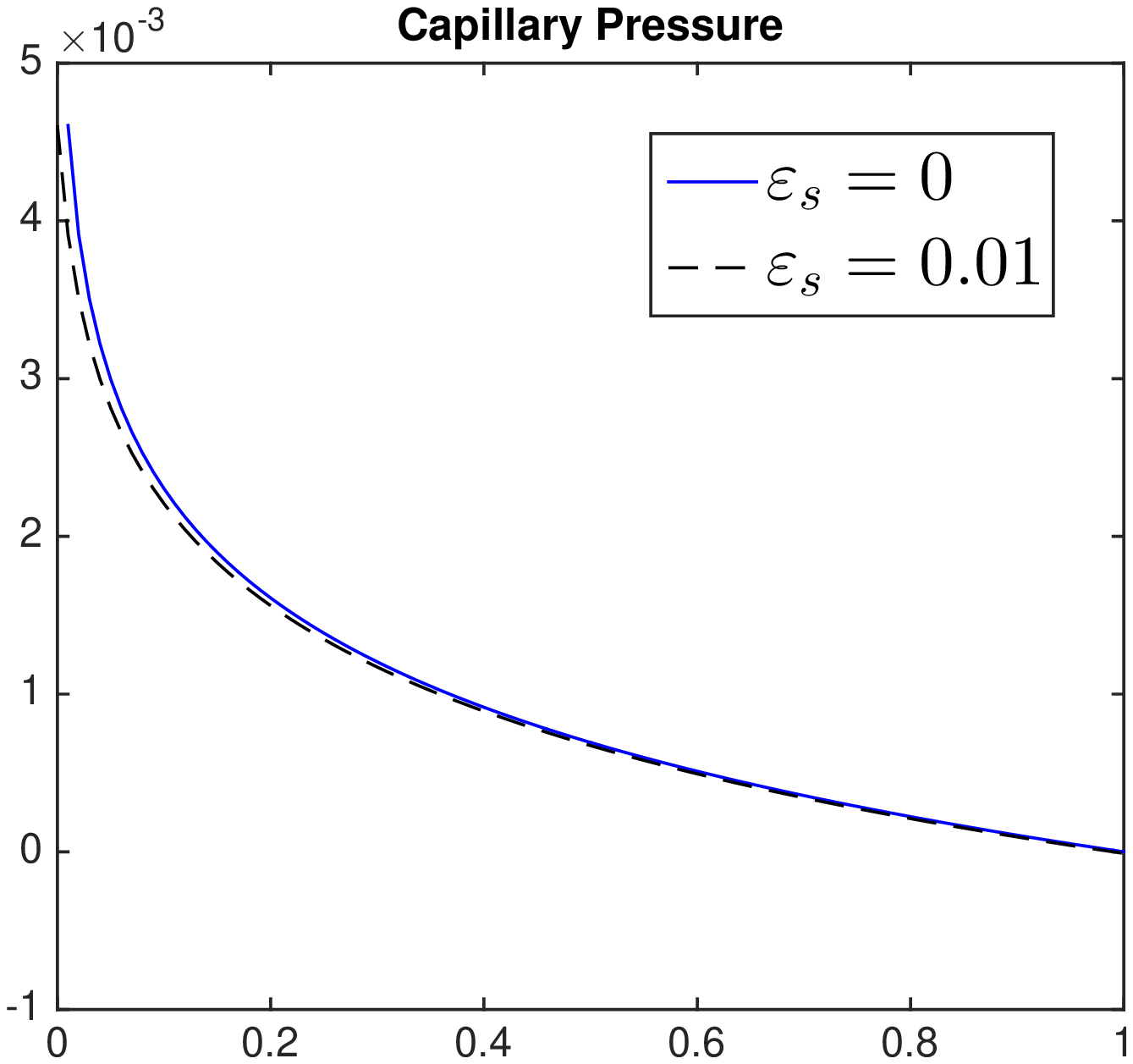}
\label{fig:num_perm_a}
}
\hspace*{0.05in}
\subfloat[Relative permeabilities]
{
\includegraphics[width=0.3\textwidth]{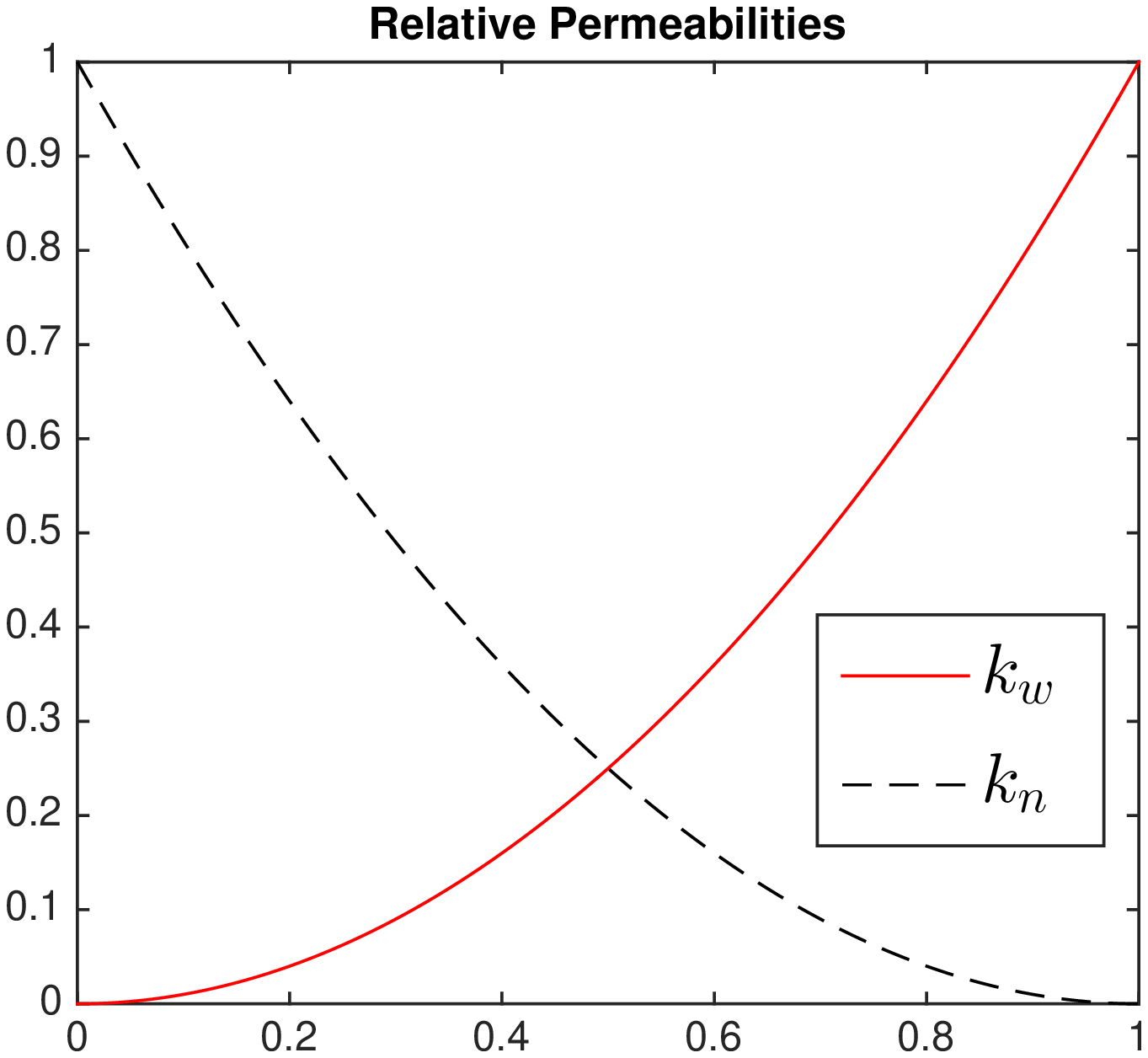}
\label{fig:num_perm_b}
}
\caption{Example 1. Given capillary pressure \eqref{eqn:capillary} values
and
relative permeabilities \eqref{ex1_relperm}.}
\label{fig:num_perm}
\end{figure}
The capillary pressure is defined as
\begin{equation}
p_c(s_w) := \dfrac{B_c}{\sqrt{K}} \log({s}_w+\varepsilon_s),
\label{eqn:capillary}
\end{equation}
where
$K$ is the absolute permeability in Darcy scale (i.e $\SI{1}{D} = \SI{9.869233e-13}{m^2})$
and
$K =  K_D I$ with $K_D=\SI{e-5}{D}$, 
where $I$ is an identity matrix, $B_c = -{0.0001}$ and 
$\varepsilon_s = 0.01$ to avoid zero singularity (see Figure \ref{fig:num_perm_a} ). 
If $s_w + \varepsilon_s \geq 1$ then we set to $s_w + \varepsilon_s=1$.
Relative permeabilities are given as a function of the wetting phase saturation, 
\begin{equation}
k_w(s_w) := s_w^2, \  \text{ and }   \ k_n(s_w) := (1-s_w)^2 ;
\label{ex1_relperm}
\end{equation}
see Figure \ref{fig:num_perm_b} for more details.
In addition, we define following the parameters:
$\mu_w = \SI{1}{cp}$, $\mu_n = \SI{2}{cp}$,
$\rho_w = \rho_n = \SI{1000}{kg/m^3}$,
$\bg=[0,\SI{-9.8}{m/s^2}]/101325$ 
(scaling with pressure (atm) $\SI{1}{atm} = \SI{101325}{\pascal}$), 
$c_w^F=\num{e-12}$, and $\phi = 0.8$.

We illustrate the convergence of EG flow \eqref{eq:egscheme} and EG saturation \eqref{eq:egscheme_saturation}, separately for the two-phase flow system with capillary pressure. 
In this case, exact values of $s_w(t^{k})$ and $s_w(t^{k-1})$ are provided to compute $P_w^{k+1}$, and exact values of $p_w(t^{k})$ and $p_w(t^{k-1})$ are provided to compute each $S_w^{k+1}$.
The entropy residual stabilization term \eqref{eqn:discrete_transport_ev} discussed in Section \ref{sec:levelset_entropy} is included with 
$\lambda_{\textsf{Ent}} =\lambda_{\textsf{Lin}}  =\num{e-2}$ 
and entropy function \eqref{eqn:entropy_func_2} chosen with $\varepsilon = \num{e-4}$.
The penalty coefficients are set as $\alpha = 100$ and $\alpha_T = 0.01$. 
For each of the flow and transport equations, respectively, 
five computations on uniform meshes were computed where the mesh size $h$ is divided by two for each cycle.
The time discretization is chosen fine enough not to influence the spatial errors and the time step $\Delta t$ is divided by two for each cycle. 
Each cycle has $100,200,400,800$ and $\num{1600}$ time steps and the errors are computed at the final time $\mathbb{T}=0.1$.
\begin{figure}[!h]
\centering
\subfloat[Pressure Error]
{
\includegraphics[width=0.35\textwidth]{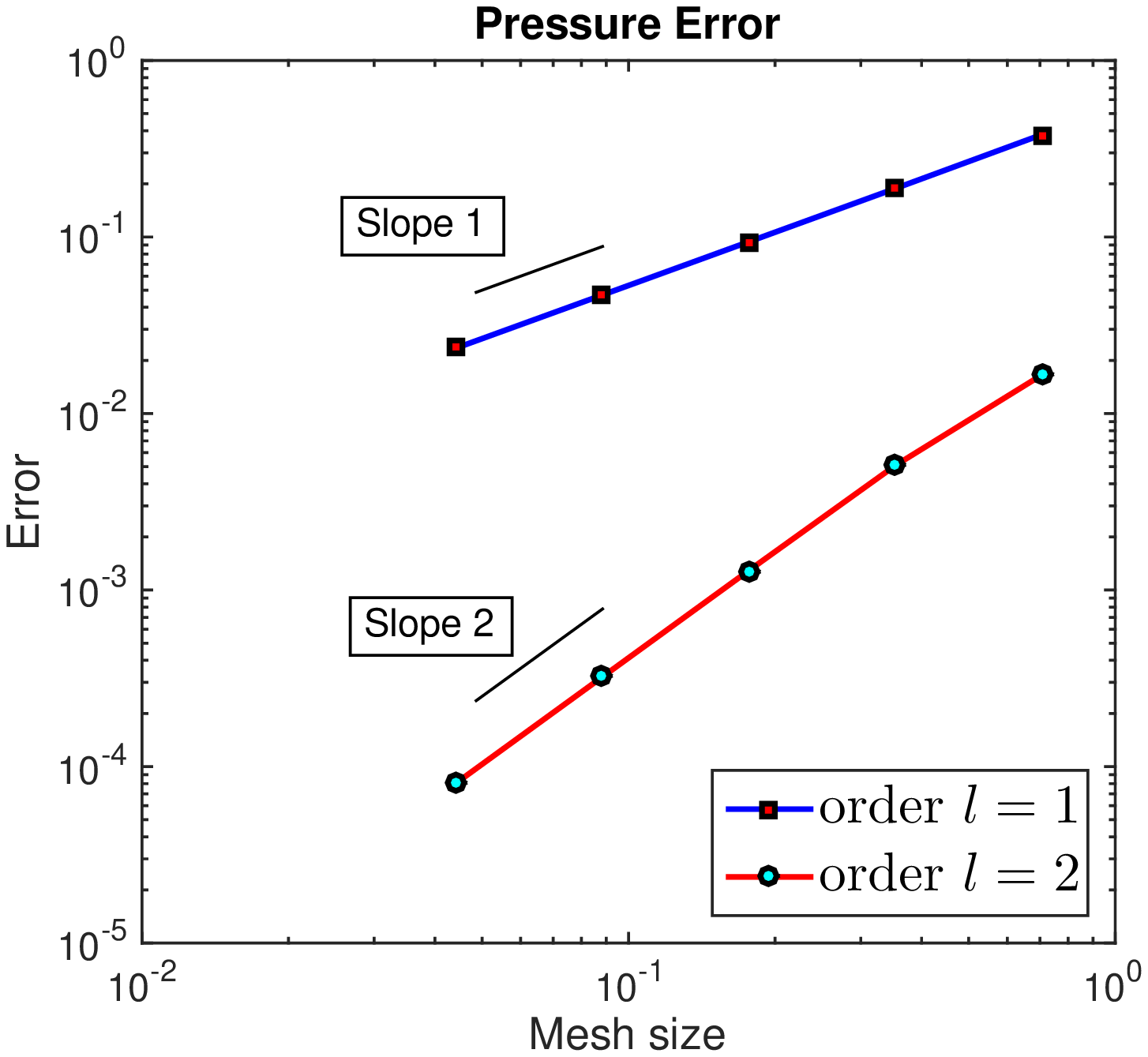}
\label{fig:ex_1}
}
\subfloat[Saturation Error]
{
\includegraphics[width=0.35\textwidth]{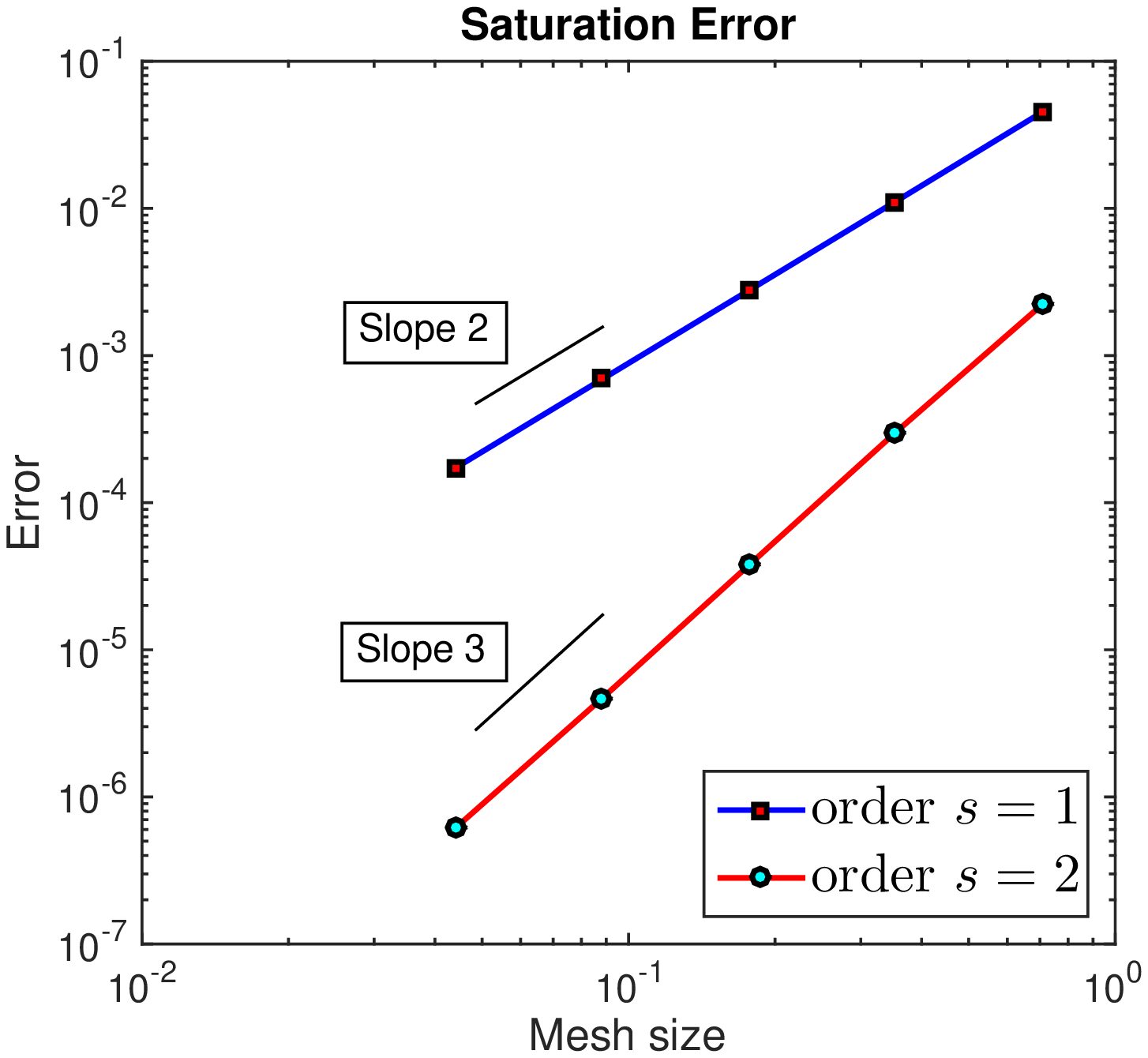}
\label{fig:ex_1_s}
}
\caption{Example 1. Decoupled case. Error convergence rates for pressure and saturation in semi-$H^1$ norm and $L^2$ norm, respectively.
Optimal order of convergences are observed for both linear and quadratic order cases.}
\end{figure}

The behavior of the $H^1(\Omega)$ semi norm errors for the approximated pressure solution versus the mesh size $h$ are depicted in Figure \ref{fig:ex_1}. 
Next, the $L^2(\Omega)$  error for the approximated saturation solutions versus the mesh size is illustrated in Figure \ref{fig:ex_1_s}. 
Both linear and quadratic orders ($l,s=1,2$) were tested and the optimal order of convergences as discussed in \cite{LeeLeeWhi15} are observed.

\subsection{Example 2. Convergence Tests - coupled case}
\begin{figure}[!h]
\centering
\subfloat[Pressure Error]
{
\includegraphics[width=0.35\textwidth]{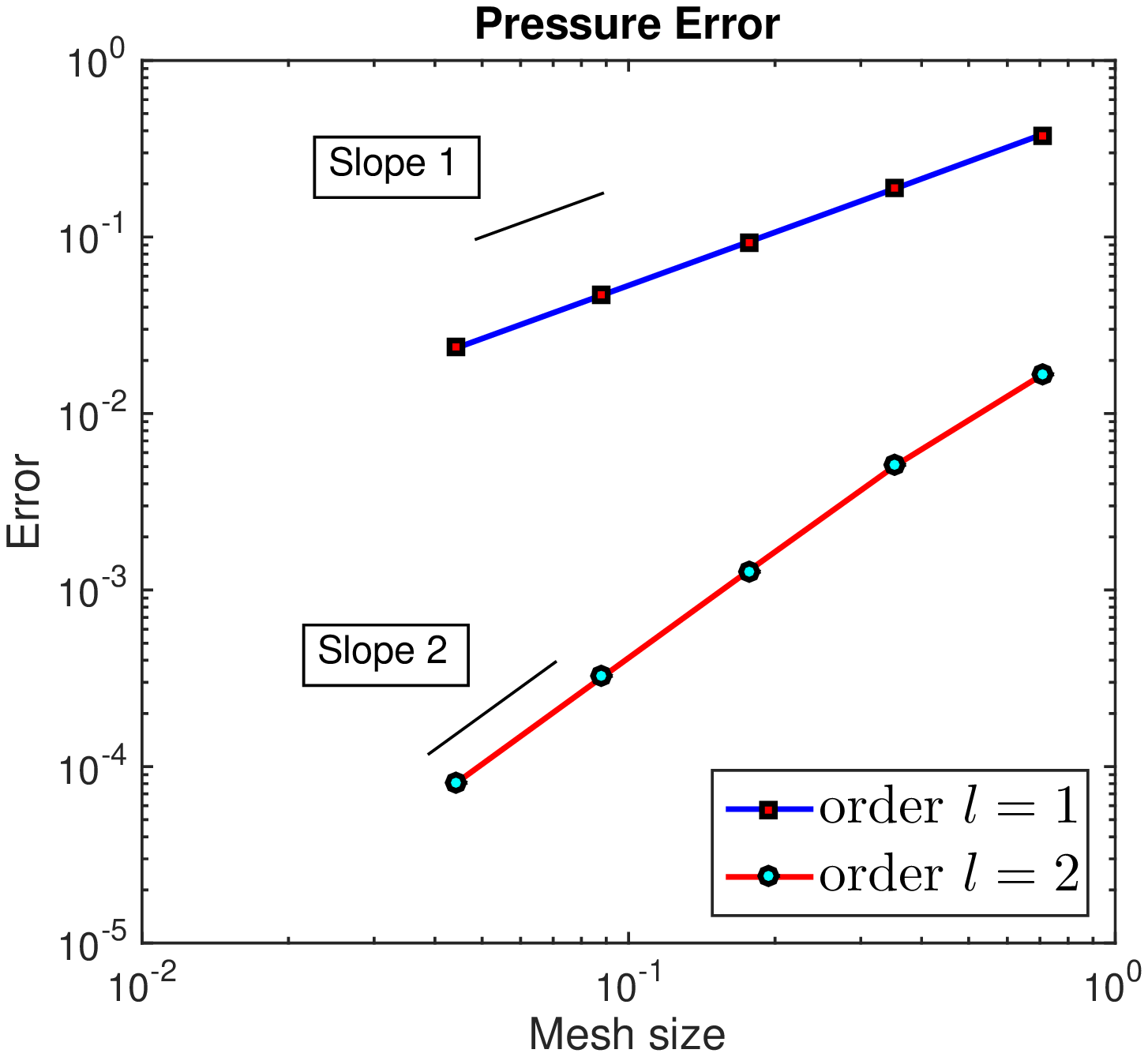}
\label{fig:ex_2}
}
\subfloat[Saturation Error]
{
\includegraphics[width=0.35\textwidth]{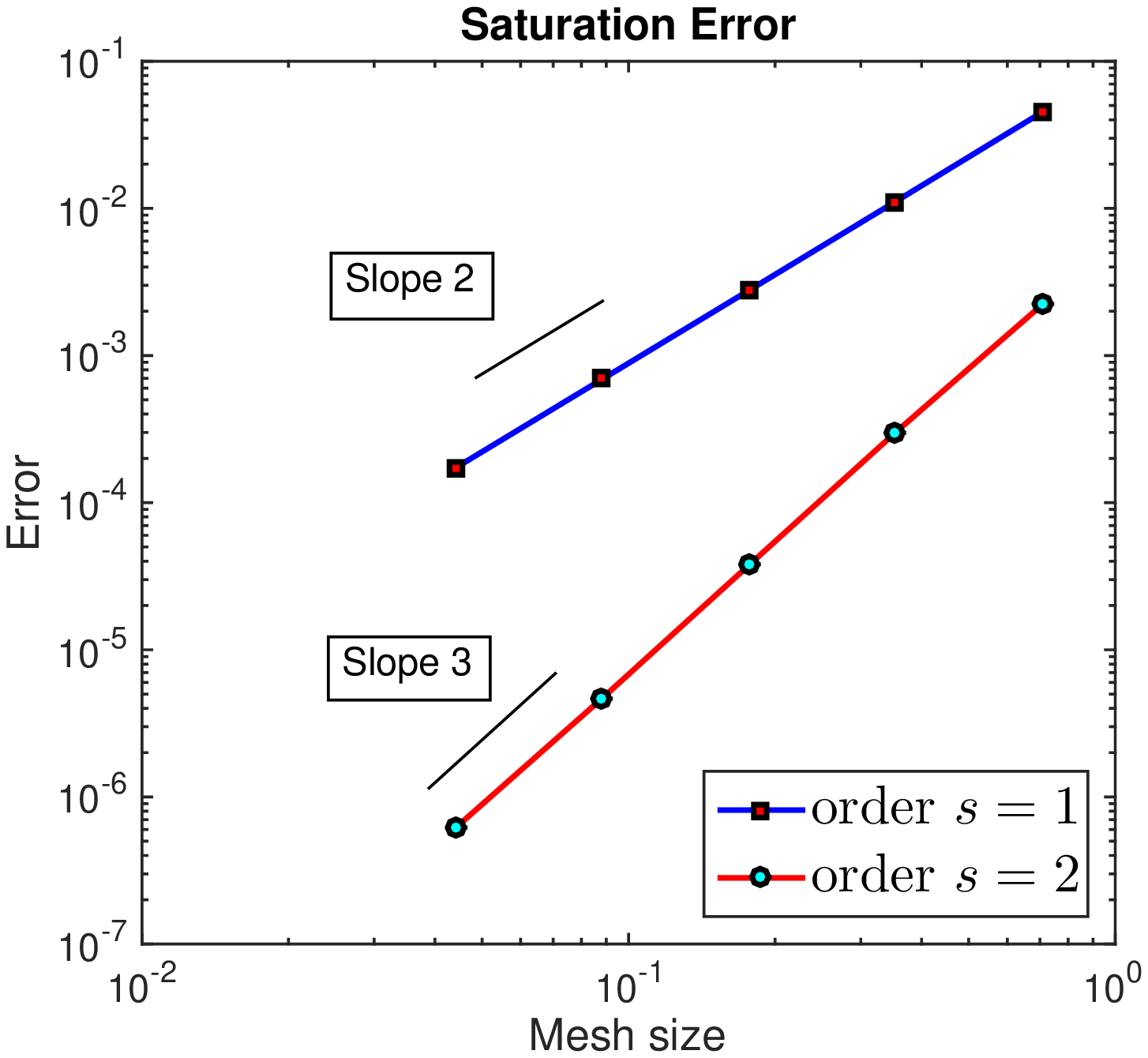}
\label{fig:ex_2_s}
}
\caption{Example 2. Coupled case (sequential IMPES). Error convergence rates for pressure and saturation in semi-$H^1$ norm and $L^2$ norm, respectively.
}
\label{fig:ex_2_total}
\end{figure}

\begin{figure}[!h]
\centering
\subfloat[Pressure Error]
{
\includegraphics[width=0.35\textwidth]{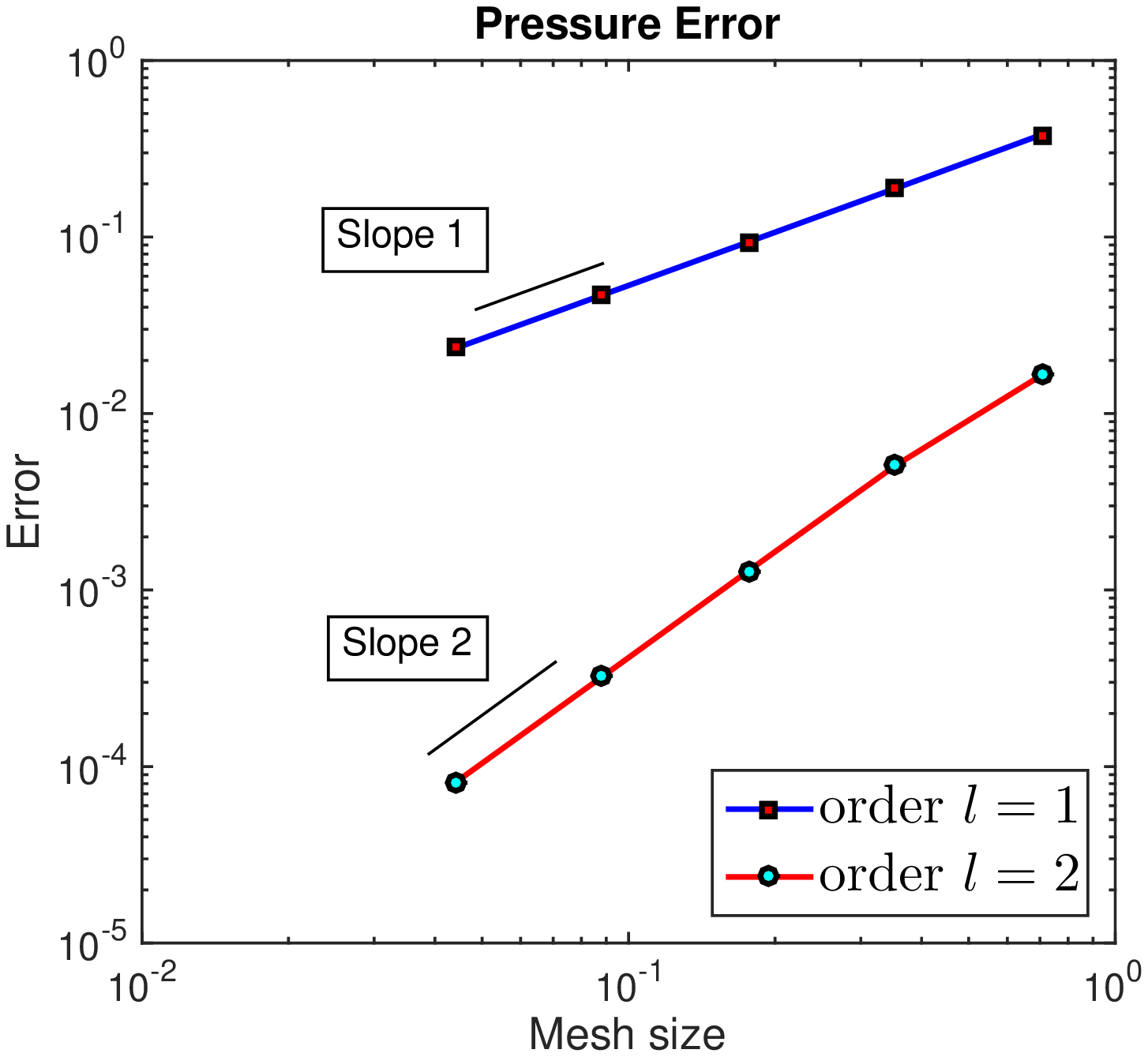}
\label{fig:ex_2-1}
}
\subfloat[Saturation Error]
{
\includegraphics[width=0.35\textwidth]{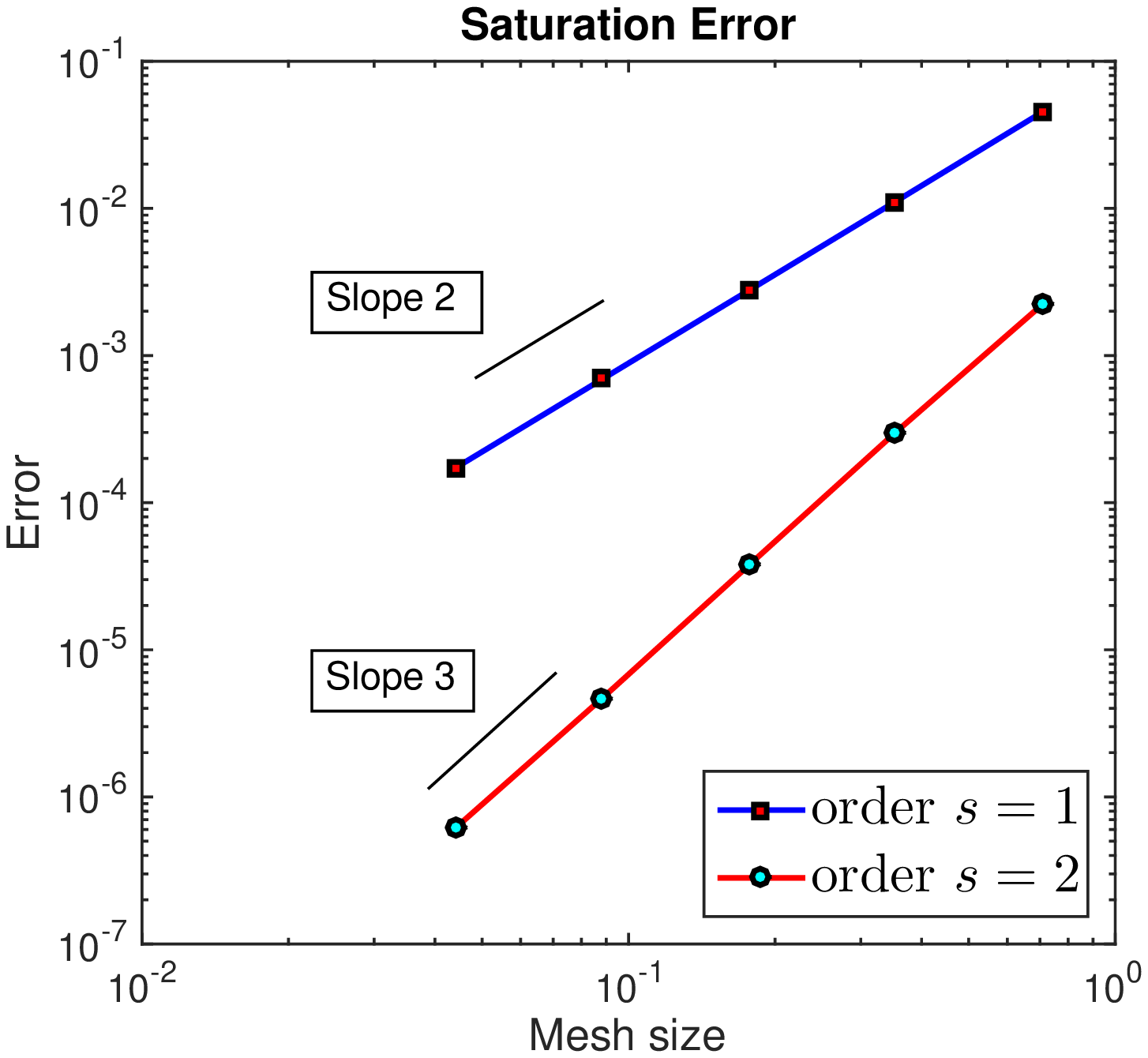}
\label{fig:ex_2-1_s}
}
\caption{Example 2. Coupled case (iterative IMPES). Error convergence rates for pressure and saturation in $H^1$ semi norm and $L^2$ norm, respectively.
}
\label{fig:ex_2-1_total}
\end{figure}

In this section, we solve the same problem as in the previous example but with a pressure and saturation system coupled.
Here, two different algorithms were tested and compared:
sequential IMPES (Section \ref{sec:IMPES}) and  iterative IMPES (Section \ref{sec:itIMPES}). 
The convergences of the errors for the pressure and the saturation are provided in Figures \ref{fig:ex_2_total} and \ref{fig:ex_2-1_total}.
We observed that the optimal rates of convergence for the high order cases  ($l=2, s=2$)
are obtained for both the sequential and  iterative IMPES scheme. 
Here the tolerance was set to $\varepsilon_{I}= \num{e-10}$ and  
3-4 iterations were required for the convergence at each time step for iterative IMPES.

\subsection{Example 3. A homogeneous channel.}
\begin{figure}[!h]
\centering
\subfloat[$s_w$ value.]
{
\includegraphics[width=0.5\textwidth]{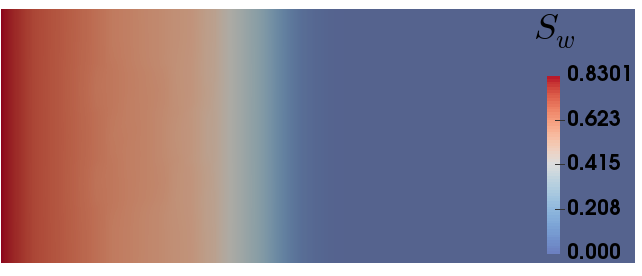}
\label{fig:ex_3_a}
}
\subfloat[$S_w$ value over line $y=0.25$. ]
{
\includegraphics[width=0.5\textwidth]{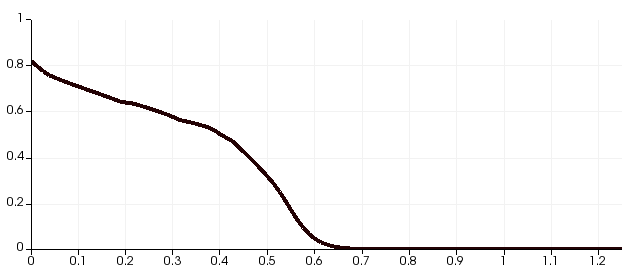}
\label{fig:ex_3_lin_ent}
}\\
\subfloat[$\mu_{{\textsf{Ent}}}(S_w, \hat{\bU}_w)_{|T}$ values.]
{
\includegraphics[width=0.5\textwidth]{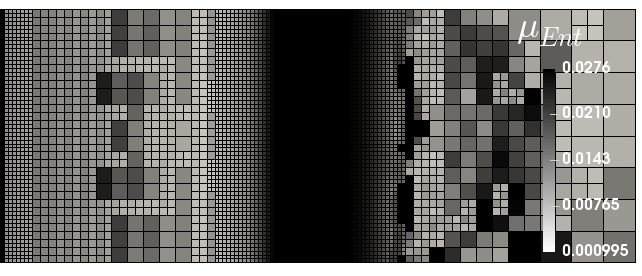}
\label{fig:ex_3_entropy}
}
\subfloat[Choices of the stabilization coefficient.]
{
\includegraphics[width=0.5\textwidth]{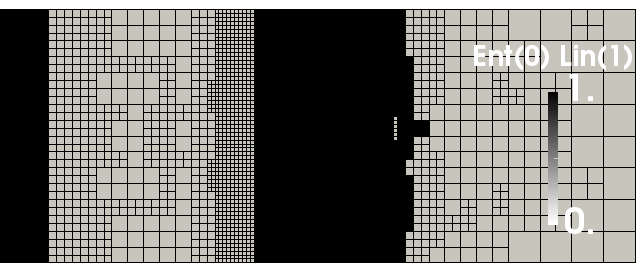}
\label{fig:ex_3_entropy_choice}
}
\caption{Example 3. Numerical results at time step number $k=50$.
(a) wetting phase saturation value. 
(b) values at (a) are plotted over the fixed line $y=0.25$.
(c) entropy residual viscosity values for each cell. 
(d) choices of viscosities; linear viscosity is chosen where the entropy residual values are larger.}
\label{fig:ex_3}
\end{figure}
In this example, we illustrate the computational features of our algorithms including entropy viscosity stabilization and dynamic mesh adaptivity with zero capillary pressure. 
The computational domain is $\Omega=({0},\SI{0}{\metre})\times(\SI{1.25}{\metre},\SI{0.5}{\metre})$ and 
the domain is saturated with a non-wetting phase, residing fluid ($s^0_n = 1$ and $s^0_w = 0$).  
A wetting phase fluid is injected at the left-hand side of the domain, thus 
$$
p_{w,\textsf{in}} = \SI{1}{atm},\  s_{w,\textsf{in}}=1 \ \text{ on } \ x=\SI{0}{\metre}.
$$ 
On the right hand side, we impose
$$
p_{w,\textsf{out}}  = \SI{0}{atm} \ \text{ on } \  x=\SI{1.25}{\metre},
$$
and no-flow boundary conditions on the top and the bottom of the domain.
Fluid and rock properties are given as 
$\mu_w = \SI{1}{cP}$, $\mu_n = \SI{3}{cP}$, 
$\rho_w = \SI{1000}{kg/m^3}$, $\rho_n=\SI{830}{kg/m^3}$, $K_D = \SI{1}{D}$,
$c_w^F = \num{e-8}$ and $\phi = 0.2$. 
Relative permeabilities are given as a function of the wetting phase saturation \eqref{ex1_relperm}, and 
the capillary pressure is set to zero for this case. 
The penalty coefficients are set as $\alpha = 100$ and $\alpha_T = 100$.

Figure \ref{fig:ex_3} illustrates the wetting phase saturation ($s_w$) at the time step number $k=50$ with the entropy stabilization coefficients  
($\lambda_{\textsf{Ent}} = 0.1$, $\lambda_{\textsf{Lin}}=1$)   
and entropy function \eqref{eqn:entropy_func_2} chosen with $\varepsilon = \num{e-4}$.
Dynamic mesh adaptivity is employed with 
initial refinement level $\textsf{Ref}_T =4$, maximum refinement level $R_{\max}=6$ and minimum refinement level  $R_{\min}=2$.
Here $C_R$ is chosen to mark and refine the cells which represent 
the top 20$\%$ of the values \eqref{eqn:ent_res} over the domain and 
$C_C$ is chosen to mark and 
coarsen the cells which represent the bottom 5$\%$ of the values \eqref{eqn:ent_res} over the domain.
The initial number of cells was approximately $2000$ and maximum cell number was approximately $6000$
with a minimum mesh size  $h_{\min}=\num{1.1e-02}$.
The uniform time step size was chosen as $\Delta t = \num{5e-3}$ (CFL constant around $0.5$).
Figure \ref{fig:ex_3_lin_ent} plots the values of $S_w$ over the fixed line $y=\SI{0.25}{\metre}$.
We observe a saturation front without any spurious oscillations. 
In addition, 
Figure \ref{fig:ex_3_entropy} presents the 
adaptive mesh refinements and entropy residual values \eqref{eqn:ent_vis} at the time step number $k=50$. 
This choice of stabilization \eqref{eqn:levelset_stab_v} 
performs as expected;
see  Figure \ref{fig:ex_3_entropy_choice}.
We note that the linear viscosity \eqref{eqn:level:fv} is chosen where the entropy residual values are larger.

\subsection{Example 4. A layered three dimensional domain}
\begin{figure}[!h]
\centering
\includegraphics[width=0.45\textwidth]{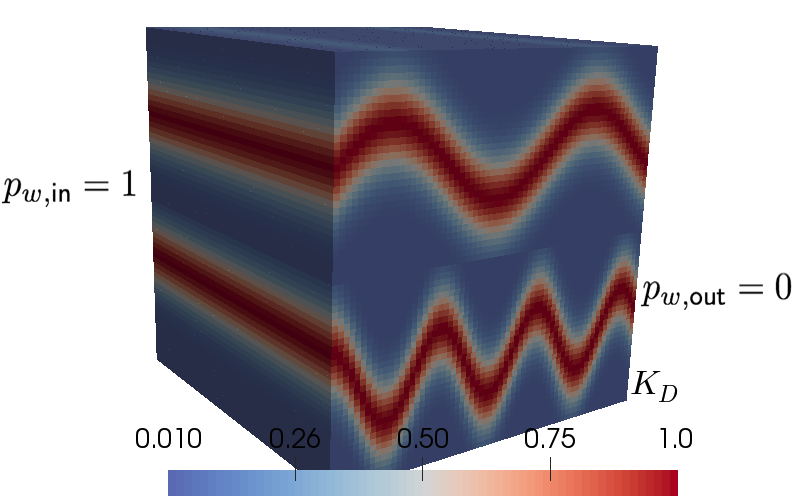}
\caption{Example 4. Setup with a given permeability ($K_D$ values).}
\label{fig:ex_4_set}
\end{figure}

This example presents a three dimensional computation in 
$\Omega=(0,\SI{1}{\metre})^3$ with a given heterogeneous domain, see Figure \ref{fig:ex_4_set} for details and boundary conditions.
Permeabilities are defined as 
$K_D= \max( \exp(-d_1^2/0.01), 0.01) $, where 
$d_1 = | y -0.75 - 0.1 * \sin (10 x) |$ for $y>0.5$ and 
$K_D= \max( \exp(-d_2^2/0.01), 0.01) $, where 
$d_2 = | y -0.25 - 0.2 * \sin (x) |$ for $y<0.5$.
All other physical parameters are the same as in the previous example.

\begin{figure}[!h]
\centering
\subfloat[$k=10$]
{
\includegraphics[width=0.33\textwidth]{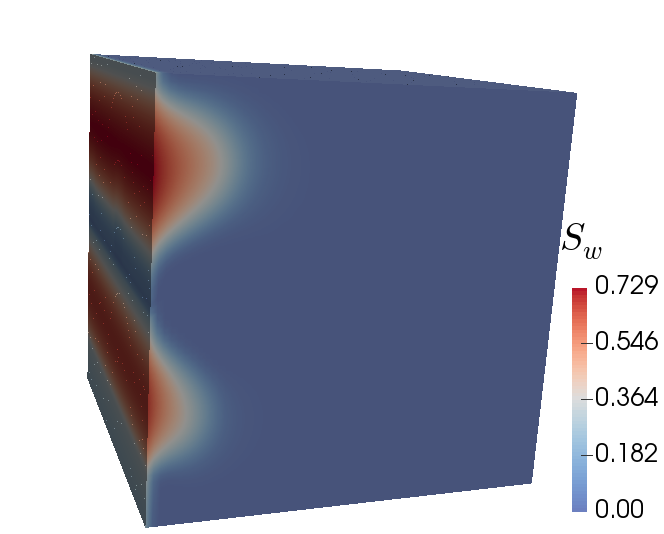}
\label{fig:ex_4_num_a}
}
\subfloat[$k=10$ with mesh]
{
\includegraphics[width=0.33\textwidth]{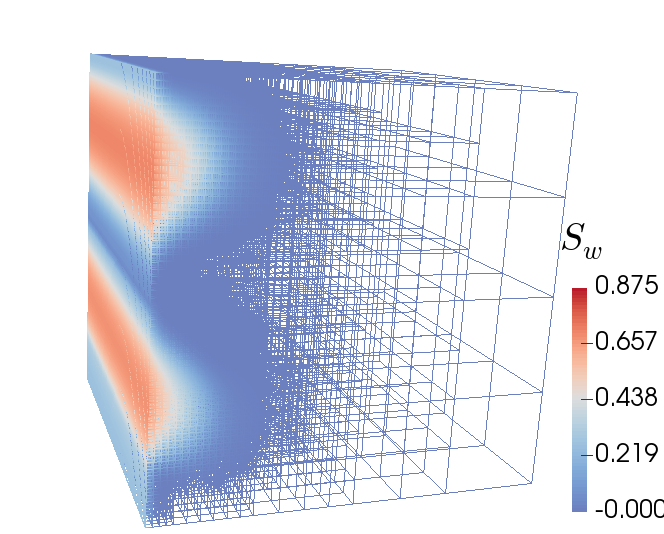}
\label{fig:ex_4_num_mesh_a}
}
\subfloat[$k=50$]
{
\includegraphics[width=0.33\textwidth]{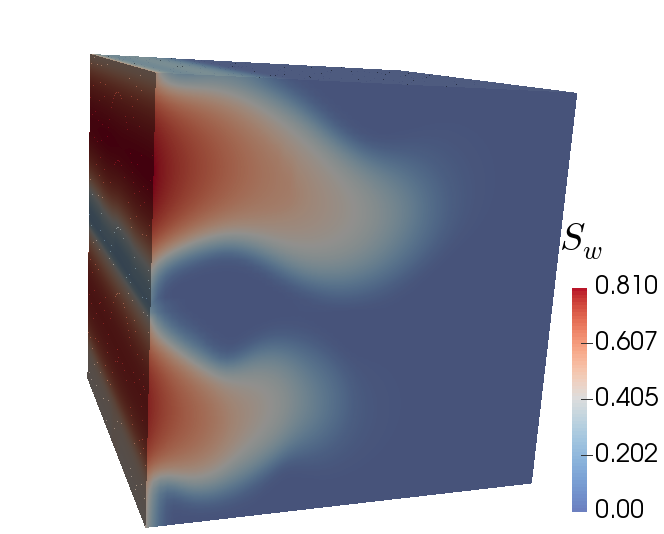}
\label{fig:ex_4_num_b}
}\\
\subfloat[$k=100$ with mesh] 
{
\includegraphics[width=0.33\textwidth]{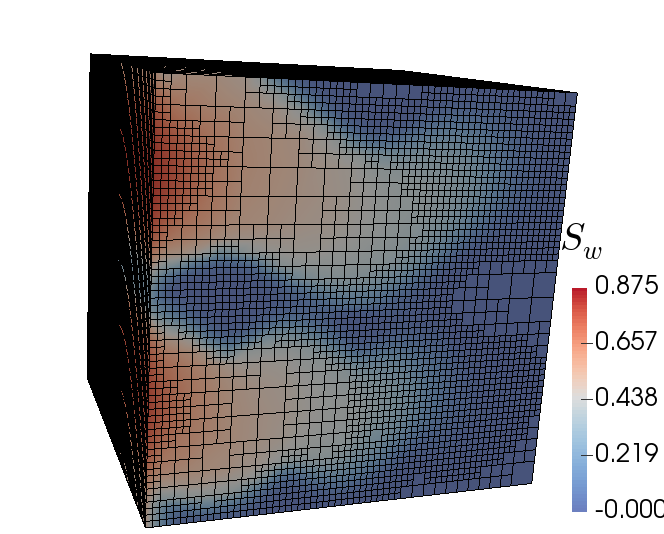}
\label{fig:ex_4_num_mesh_b}
}
\subfloat[$k=200$]
{
\includegraphics[width=0.33\textwidth]{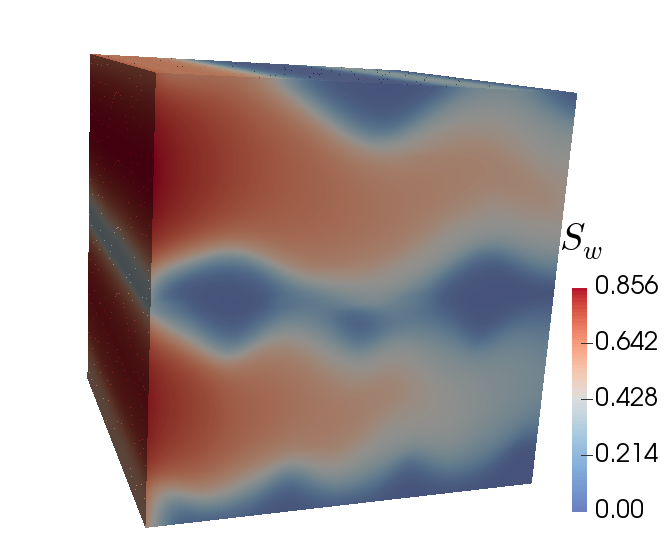}
\label{fig:ex_4_num_c}
}
\subfloat[$k=300$]
{
\includegraphics[width=0.33\textwidth]{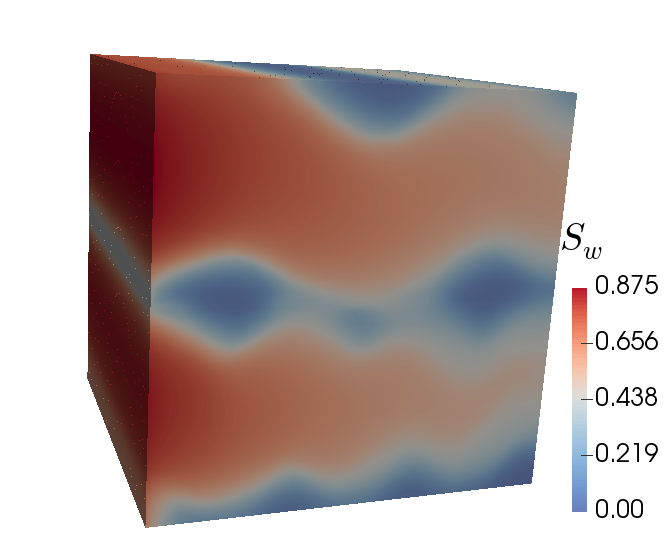}
\label{fig:ex_4_num_d}
}
\caption{Example 4. 
The wetting phase saturation ($S_w$) at each time step number with adaptive mesh refinements. }
\label{fig:ex_4_num}
\end{figure}

Figure \ref{fig:ex_4_num}  illustrates the wetting phase saturation ($S_w$) at the time step number $k=10, 50, 200, $ and $300$ with the entropy stabilization coefficients  
($\lambda_{\textsf{Ent}} = 0.25$, $\lambda_{\textsf{Lin}}=0.5$)  
and entropy function \eqref{eqn:entropy_func_2} chosen with $\varepsilon = \num{e-3}$.
Dynamic mesh adaptivity is employed with 
$\textsf{Ref}_T =4$,  $R_{\max}=6$ and  $R_{\min}=2$.
The number of cells at $k=300$ is around $\num{262100}$ and 
the minimum mesh size is $h_{\min}=\num{0.027}$ with a time step size $\Delta t = 0.006$ (CFL constant is 1). 
See figures \ref{fig:ex_4_num_mesh_a}-\ref{fig:ex_4_num_mesh_b} for adaptive mesh refinements for different time steps.
The adaptive mesh refinement strategy becomes very efficient for large-scale three dimensional problems using parallelization. 

\subsection{Example 5. A benchmark: effects of capillary pressure} 
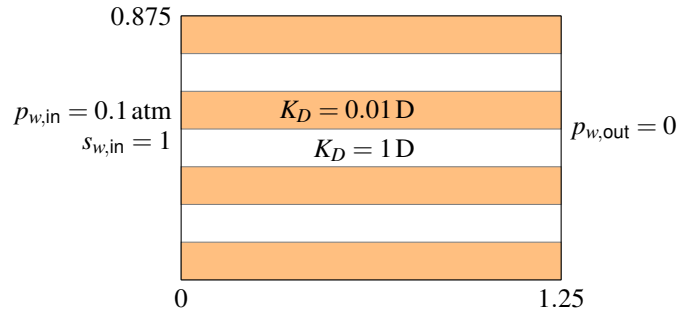
\begin{figure}[!h]
\centering
\begin{tikzpicture}[scale=4.]
\draw (0,0) -- (1.25,0); 
\draw (1.25,0) -- (1.25,0.875); 
\draw (1.25,0.875) -- (0,0.875); 
\draw (0,0.875) -- (0,0);

\node[left] at (0.,0.55)    {$p_{w,\textsf{in}}=\SI{0.1}{atm}$}; 
\node[left] at (0.,0.45)    {$s_{w,\textsf{in}}=1$}; 
\node[right] at (1.25,0.5) {$p_{w,\textsf{out}}=0$} ; 

\node[left] at (0.,0.875) {$0.875$} ; 
\node[below] at (0.,0) {$0$} ; 
\node[below] at (1.25,0) {$1.25$} ; 

\draw [fill=orange, opacity=0.5] (0.,0.75) rectangle (1.25,0.875);
\draw [fill=orange,opacity=0.5] (0.,0.5) rectangle (1.25,0.625);
\draw [fill=orange,opacity=0.5] (0.,0.25) rectangle (1.25,0.375);
\draw [fill=orange,opacity=0.5] (0.,0.) rectangle (1.25,0.125);

\node[left] at (0.8,0.56)    {$K_D=\SI{0.01}{D}$}; 
\node[left] at (0.8,0.43)    {$K_D=\SI{1}{D}$}; 
\end{tikzpicture}
\caption{Example 5. Two dimensional domain with heterogeneous permeabilities. Layered setup to test the effect of the capillary pressure. Permeabilities are defined as $K_D=\SI{0.01}{D}$ for the dark region and  $K_D=\SI{1}{D}$ for the white region.}
\label{fig:ex2_domain}
\end{figure}

In this example, we emphasize the effects of capillary pressure in a heterogeneous media as shown in \cite{Hoteit:2008dt,yang2017nonlinearly}. Here, we impose layers of different permeabilities in the 
computational domain $\Omega=(\SI{0}{\metre},\SI{0}{\metre})\times(\SI{1.25}{\metre},\SI{0.875}{\metre})$. 
See Figure \ref{fig:ex2_domain}.
\begin{figure}[!h]
\centering
\subfloat[$k=125$]
{
\includegraphics[width=0.33\textwidth]{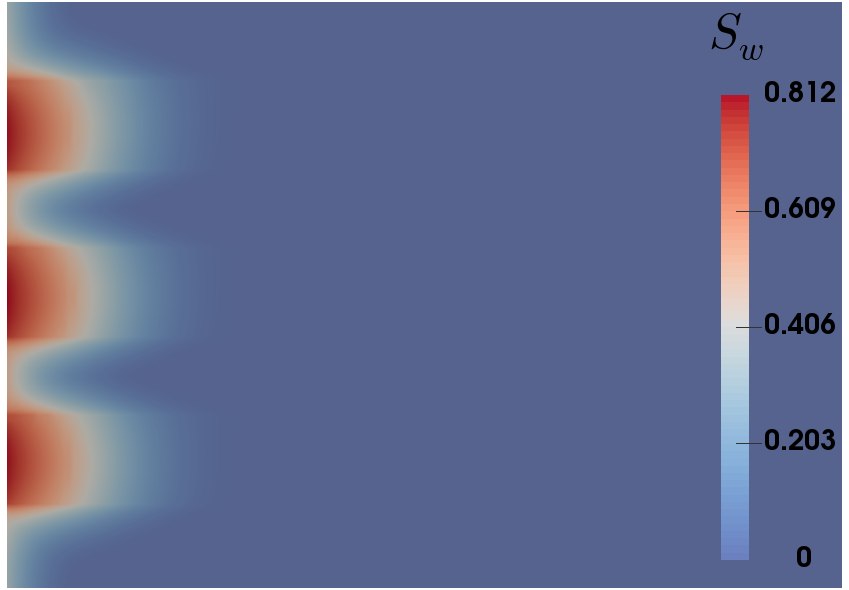}
}
\hspace*{0.01in}
\subfloat[$k=125$]
{
\includegraphics[width=0.33\textwidth]{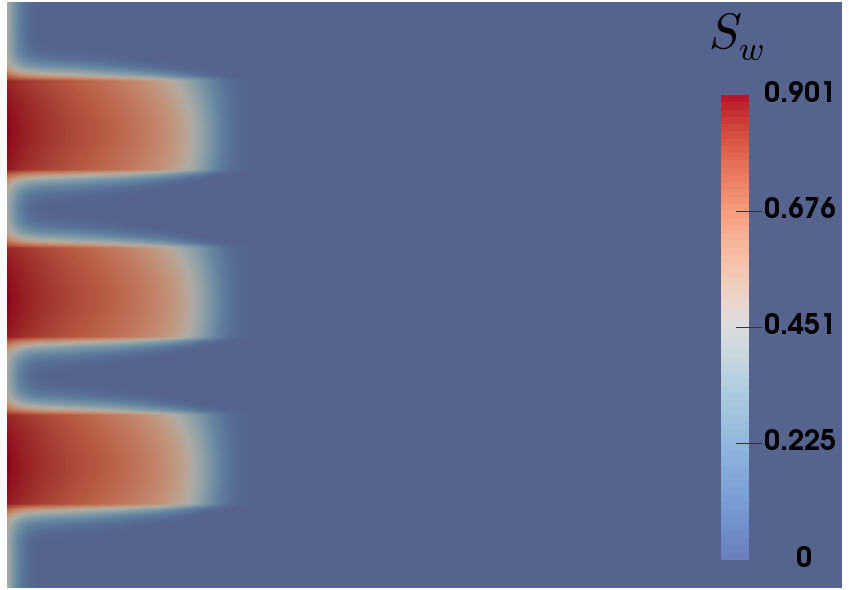}
}\\
\subfloat[$k=375$]
{
\includegraphics[width=0.33\textwidth]{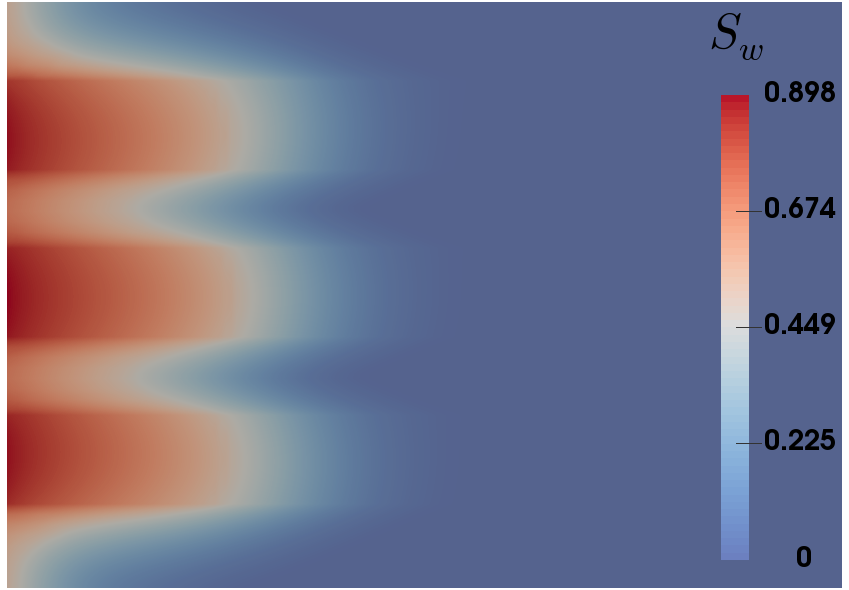}
}
\hspace*{0.01in}
\subfloat[$k=375$]
{
\includegraphics[width=0.33\textwidth]{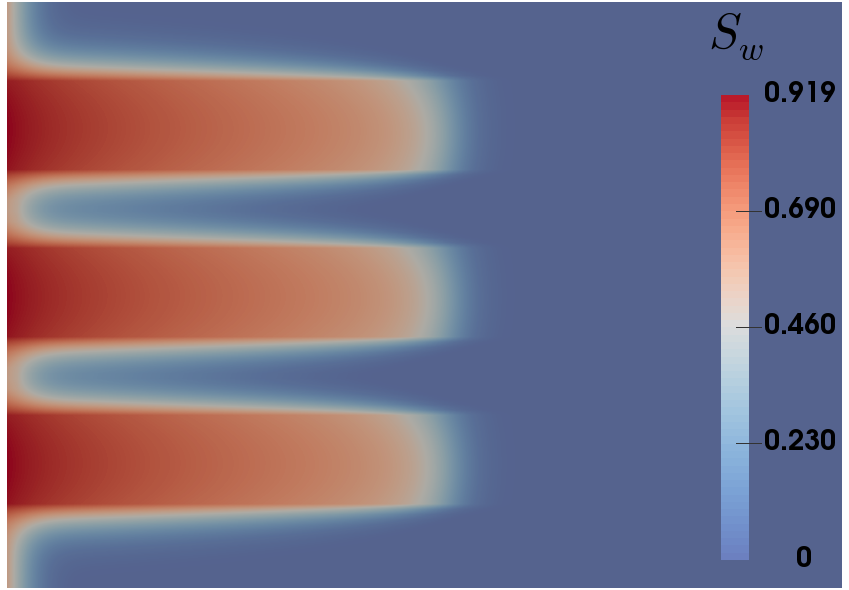}
}\\
\subfloat[$k=625$]
{
\includegraphics[width=0.33\textwidth]{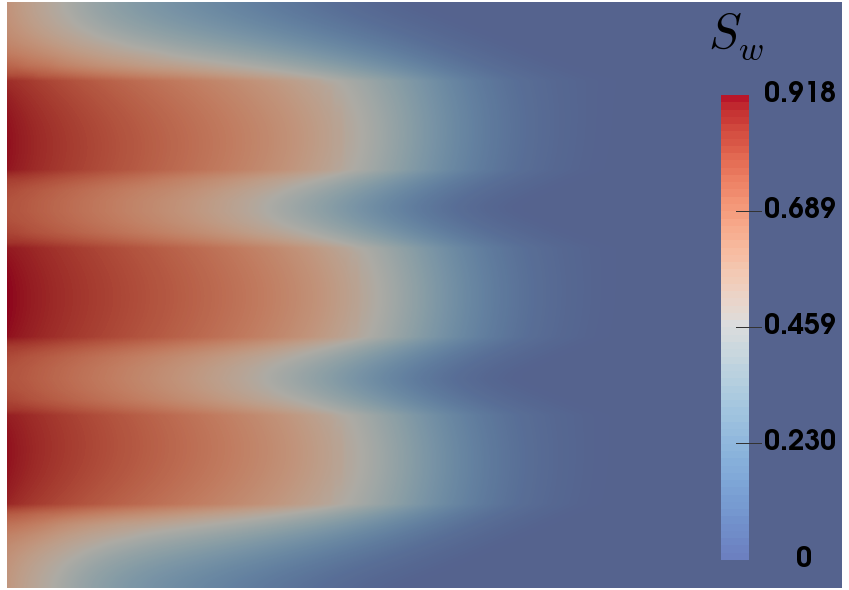}
}
\hspace*{0.01in}
\subfloat[$k=625$]
{
\includegraphics[width=0.33\textwidth]{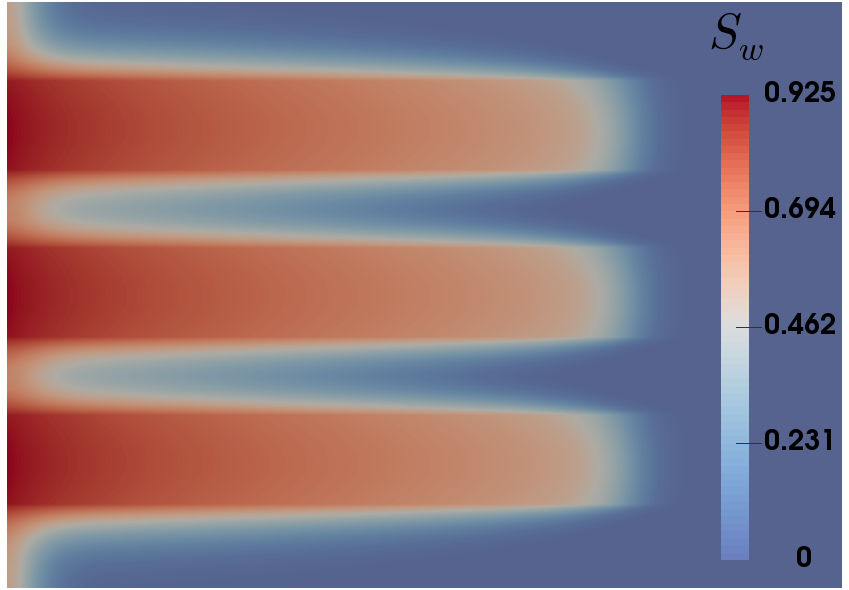}
}
\caption{Example 5. Wetting phase saturation values at each time step number. The left column (a),(c), and (e) are the values with the capillary pressure and the right column (b),(d), and (f) are the values without the capillary pressure.}
\label{fig:ex_4}
\end{figure}
The domain is saturated with a non-wetting phase (oil), 
i.e $s^0_n = 1$ and $s^0_w = 0$.
A wetting phase fluid is injected at the left-hand side of the domain, thus 
$$
p_{w,\textsf{in}} = \SI{0.1}{atm},\  s_{w,\textsf{in}}=1 \ \text{ on } \ x=\SI{0}{\metre}.
$$ 
On the right hand side, we impose
$$
p_{w,\textsf{out}}  = \SI{0}{atm} \ \text{ on } \  x=\SI{1.25}{\metre},
$$
and no-flow boundary conditions on the top and the bottom of the domain.
Fluid properties are set as 
$\mu_w = \SI{1}{cP}$, $\mu_n = \SI{0.45}{cP}$, 
$\rho_w = \SI{1000}{kg/m^3}$, $\rho_n=\SI{660}{kg/m^3}$, 
$c_w^F = \num{e-8}$, $\phi = 0.2$, and 
{$K_D = \SI{1}{D}$ or $K_D = \SI{0.01}{D}$} as illustrated in Figure \ref{fig:ex2_domain}.
Relative permeabilities are given as a function of the wetting phase saturation \eqref{ex1_relperm}, and 
the penalty coefficients are set as $\alpha = 1$, $\alpha_c = 1$ and $\alpha_T = 1000$. 
The entropy stabilization coefficients are 
$\lambda_{\textsf{Ent}} = 1$ and $\lambda_{\textsf{Lin}}=1$. 
Dynamic mesh adaptivity is employed as same as the example 3 and 
the minimum mesh size is 
$h_{\min}=\num{0.0027}$. 
The uniform time step size is taken as $\Delta t = \num{0.005}$.
The capillary pressure \eqref{eqn:capillary} is given with {$B_c = -0.01$ and $\varepsilon_s=0.1$}.

Here two tests are performed, one with  
the capillary pressure ($B_c = -0.01$) and 
a second with zero capillary pressure ($B_c = 0$). 
The differences and effects of capillary pressure are depicted at Figure \ref{fig:ex_4} for  different time steps.
The injected wetting phase water flows faster in the high permeability layers but is more diffused
in the case with capillary pressure as shown in previous results \cite{Hoteit:2008dt,yang2017nonlinearly}.
One can observe the capillary pressure is a non-linear diffusion source term for the residing non-wetting phase. This causes more uniformed movement of the injected fluid.

\subsection{Example 6. A random heterogeneous domain with different relative permeability}
\begin{figure}[!h]
\centering
\subfloat[Permeability $K$ values]
{
\includegraphics[width=0.27\textwidth]{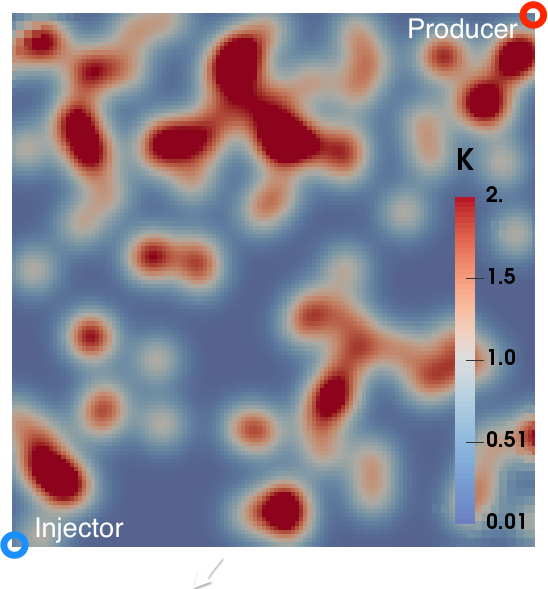}
\label{fig:ex_7_num_a}
}
\subfloat[Relative permeabilities]
{
\includegraphics[width=0.35\textwidth]{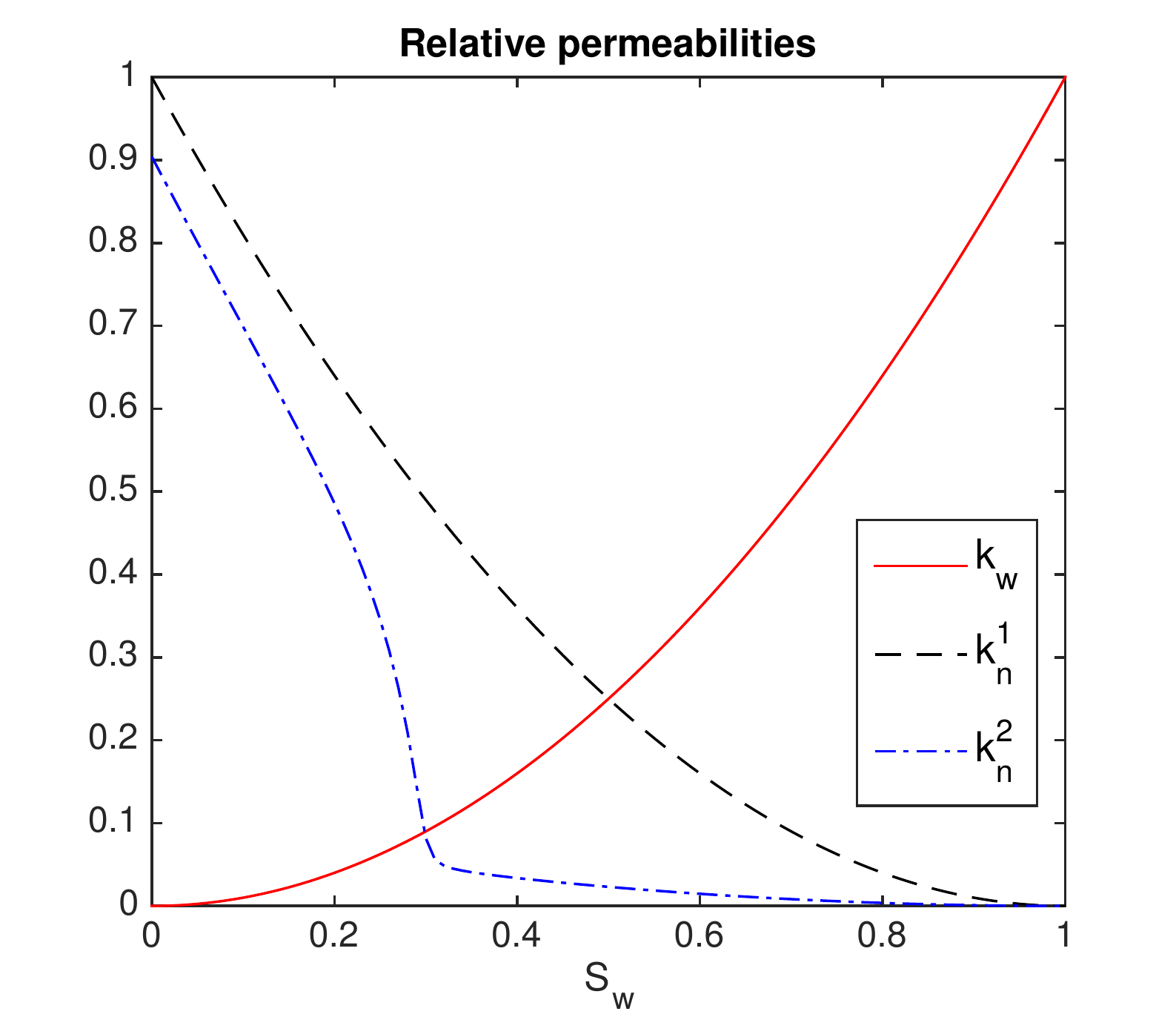}
\label{fig:ex_7_num_b} 
}
\caption{Example 6. Setup with a random absolute permeabilities, 
wetting phase relative permeability $(k_w)$, and 
two different non-wetting phase relative permeabilities ($k^1_n, k^2_n$). We note $k^2_n$ represents rough relative permeability which often  referred as the case with foam in a porous media \cite{van2017influence}.}
\label{fig:ex_7_set}
\end{figure}

This example considers well injection and production in a random heterogeneous domain 
$\Omega = (\SI{0},\SI{1}{\metre})^2$.
Wells are specified at the corners with 
injection at $(0,0)$ and production at $(\SI{1}{\metre},\SI{1}{\metre})$.
See Figure \ref{fig:ex_7_num_a} for the setup. 
We test and compare two different
non-wetting phase relative permeabilities such as 
\begin{equation}
\text{i) } 
k^1_n(s_w) := {(1-s_w)^2}
\  \text{ and }   \ 
\text{ii) }
k^2_n(s_w) := \dfrac{(1-s_w)^2}{f_w}  ,
\label{ex2_relperm}
\end{equation}
where the latter is often referred as the case with foam in a porous media \cite{van2017influence}. 
Here, 
$f_w := 1 + R (0.5 + \dfrac{1}{\pi}\arctan(\kappa (s_w - s_w^*)))$ is a mobility reduction factor  
with a constant positive parameters set to $R=10$, $\kappa = 100$, and 
a limiting water saturation $S_w^*=0.3$. 
Figure \ref{fig:ex_7_num_b} illustrates two different non-wetting phase relative permeabilities ($k^1_n, k^2_n$). The wetting phase relative permeability ($k_w$) is identical with the previous examples. 

\begin{figure}[!h]
\centering
\subfloat[$t=0.76s$ ]
{
\includegraphics[width=0.225\textwidth]{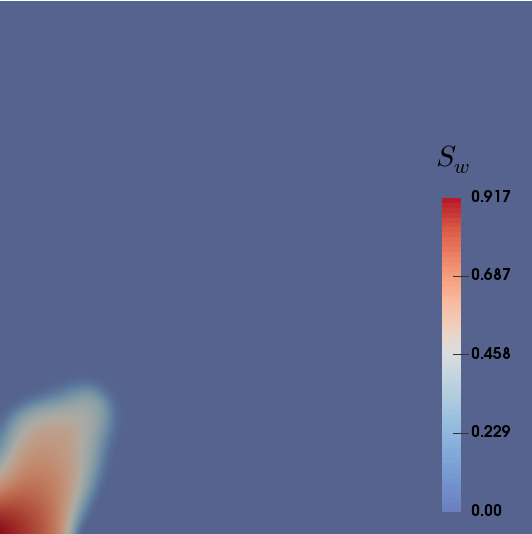}
\label{fig:ex_6_num_a}
}
\subfloat[$t=1.9s$]
{
\includegraphics[width=0.225\textwidth]{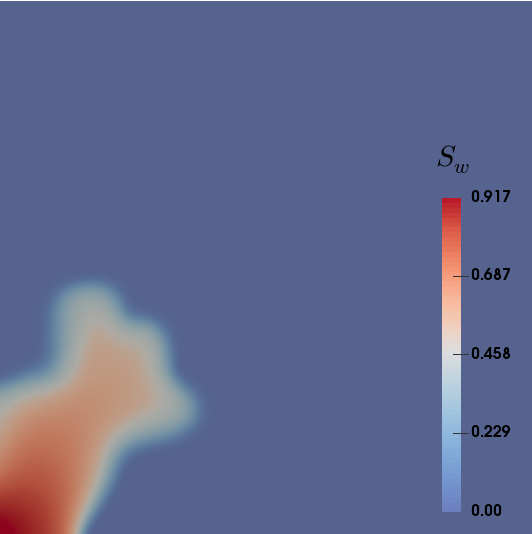}
\label{fig:ex_6_num_b} 
}
\subfloat[$t=3.04s$ ]
{
\includegraphics[width=0.225\textwidth]{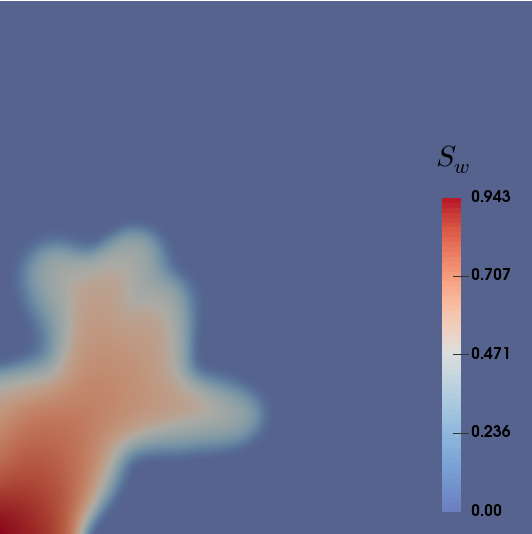}
\label{fig:ex_6_num_c}
}\\
\subfloat[$t=5s$ ]
{
\includegraphics[width=0.225\textwidth]{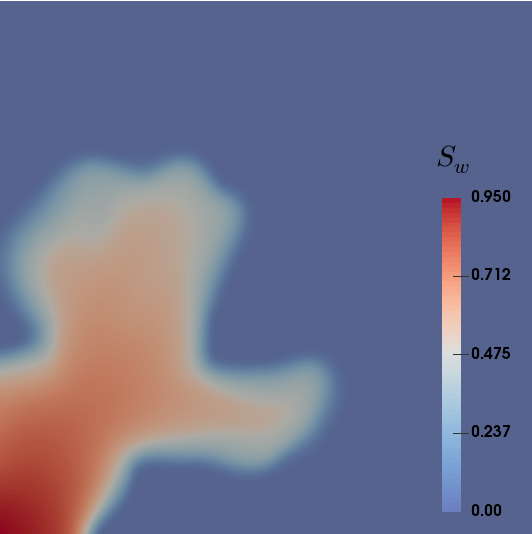}
\label{fig:ex_6_num_d}
}
\subfloat[$t=6.84s$ ]
{
\includegraphics[width=0.225\textwidth]{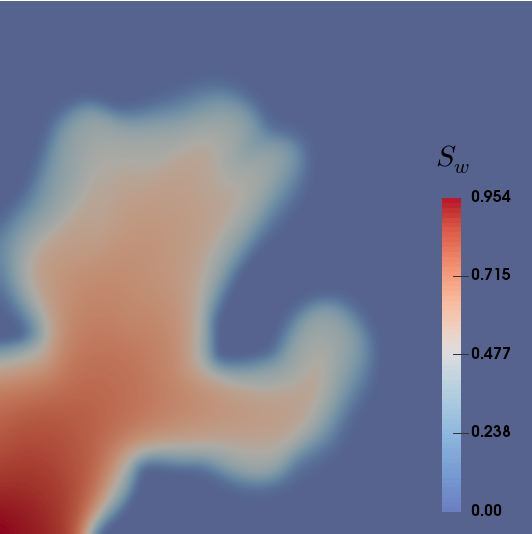}
\label{fig:ex_6_num_e}
}
\subfloat[$t=10.64s$ ]
{
\includegraphics[width=0.225\textwidth]{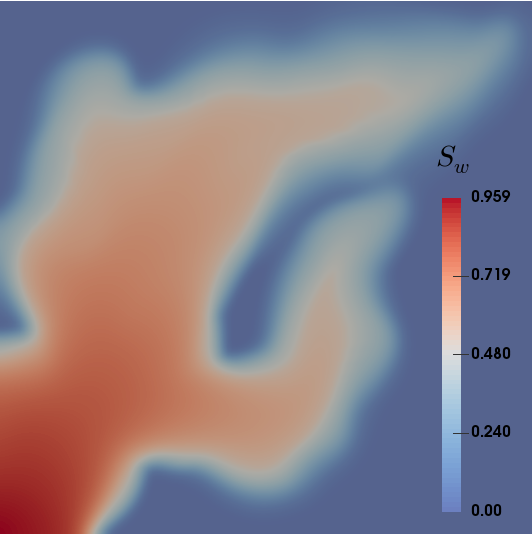}
\label{fig:ex_6_num_f}
}
\caption{Example 6. $S_w$ values for each time in a heterogeneous media with  a non-wetting phase relative permeability $k^1_n(s_w)$. 
}
\label{fig:ex_6_num}
\end{figure}
\begin{figure}[!h]
\centering
\subfloat[$t=0.765$]
{
\includegraphics[width=0.225\textwidth]{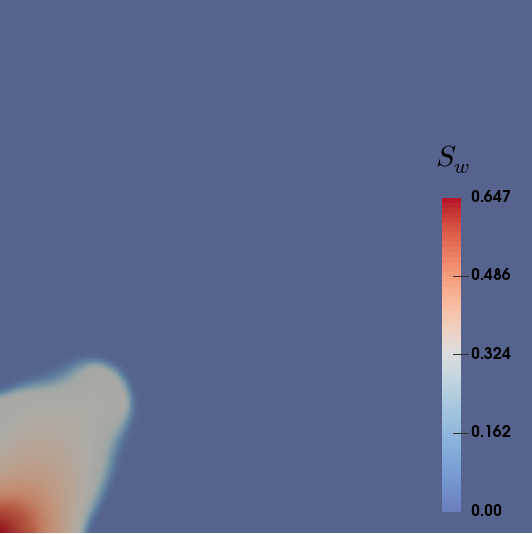}
}
\subfloat[$t=2s$]
{
\includegraphics[width=0.225\textwidth]{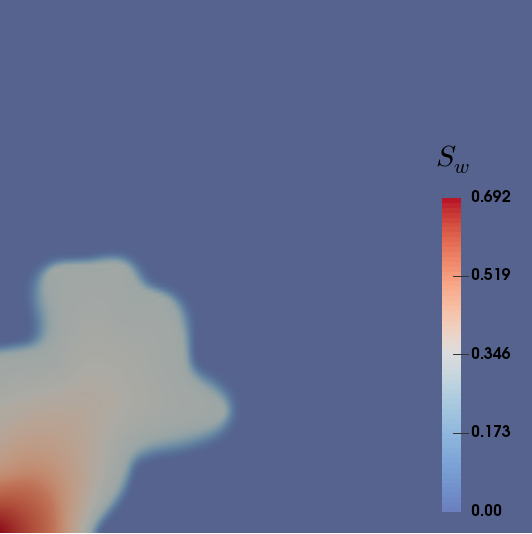}
}
\subfloat[$t=3s$]
{
\includegraphics[width=0.225\textwidth]{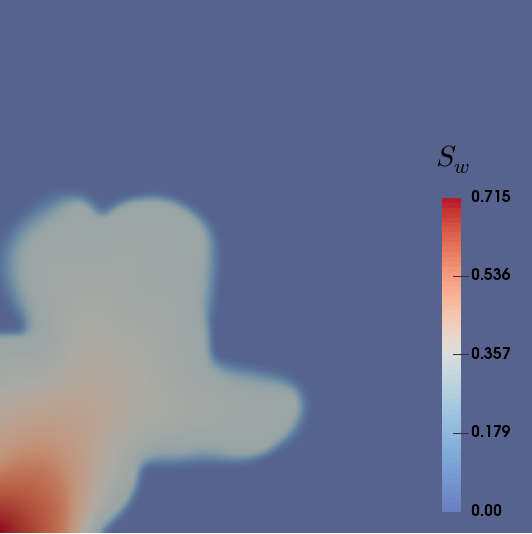}
}\\
\subfloat[$t=5s$]
{
\includegraphics[width=0.225\textwidth]{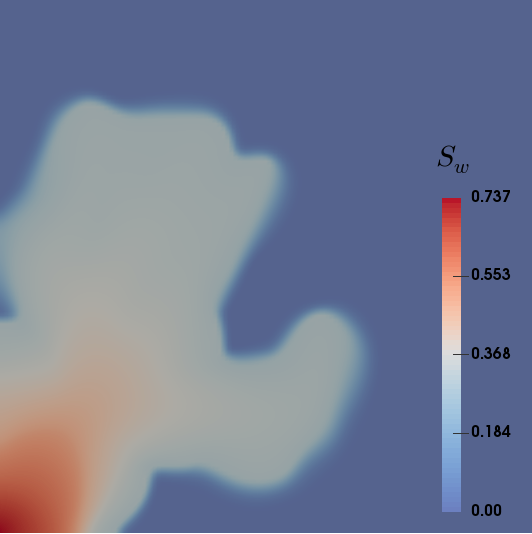}
}
\subfloat[$t=6.8s$ ]
{
\includegraphics[width=0.225\textwidth]{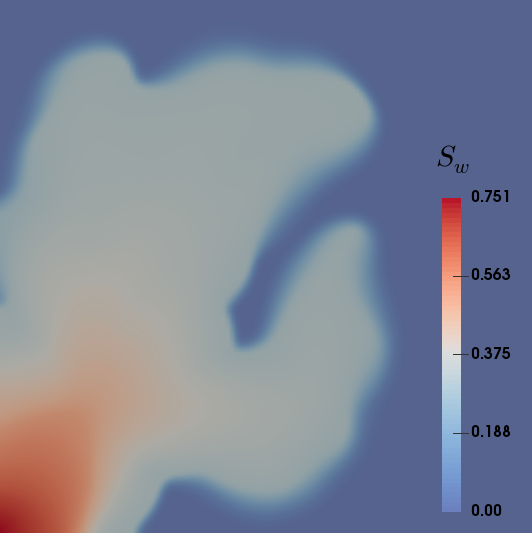}
}
\subfloat[$t=10s$ ]
{
\includegraphics[width=0.225\textwidth]{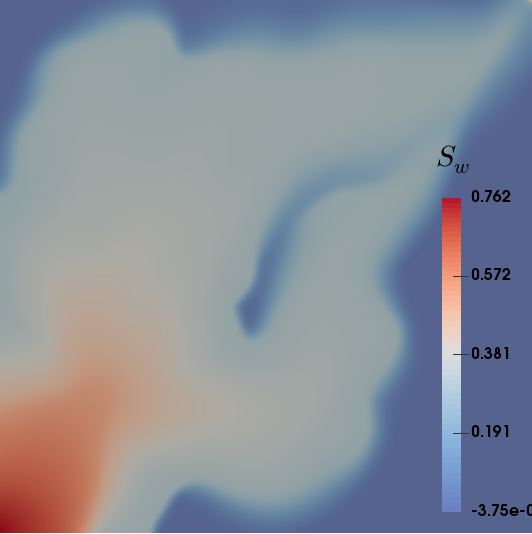}
}
\caption{Example 6. $S_w$ values for each time in a heterogeneous media with a non-wetting phase relative permeability $k^2_n(s_w)$.
}
\label{fig:ex_6_num_2}
\end{figure}

We assume the domain is saturated with a non-wetting phase, 
i.e $s^0_n = 1$ and $s^0_w = 0$ and  a wetting phase fluid  is injected.
Fluid and rock properties are given as 
$\mu_w = \SI{1}{cP}$, $\mu_n = \SI{3}{cP}$, 
$\rho_w = \SI{1000}{kg/m^3}$, $\rho_n=\SI{830}{kg/m^3}$, $c_w^F = \num{e-10}$, 
$f_w^+ = \SI{100}{m/s}$, 
$f_w^-= -\SI{100}{m/s}$, 
$f_n = 0$, 
and $\phi = 0.2$.
The capillary pressure and the gravity is neglected to emphasize the effects of heterogeneity and different non-wetting phase relative permeability. 
Here the numerical parameters are chosen as 
$h_{\min}=\num{1.1e-02}$ and $\Delta t = \num{3.8e-03}$. 
Due to the dynamic mesh refinement ($\textsf{R}_{\max} = 7$ and $\textsf{R}_{\min} = 2$), 
the number of degrees of freedom for EG transport 
and the maximum number of cells
are $32158$, $15934$, respectively at the final time $\mathbb{T}=15$. 
The entropy stabilization coefficients  are set to 
$\lambda_{\textsf{Ent}} = 0.1$ and 
$\lambda_{\textsf{Lin}}=0.25$,  
where the entropy function \eqref{eqn:entropy_func_2} is 
chosen with $\varepsilon = \num{e-3}$.
The penalty coefficients are set as $\alpha = 1$ and $\alpha_T = 1000$.

Figure \ref{fig:ex_6_num} illustrates the EG-$\polQ_1$ solution of $S_w$ values for each time in a heterogeneous media with a non-wetting phase relative permeability $k^1_n(s_w)$. Next, Figure \ref{fig:ex_6_num_2} is the case with $k^2_n(s_w)$. 
We note that wetting phase saturation values above $S_w^*$ are restricted for the latter case due to the relative permeability, $k^2_n(s_w)$.

\subsection{Example 7. A three dimensional random heterogeneous domain}
\begin{figure}[!h]
\centering
\includegraphics[width=0.35\textwidth]{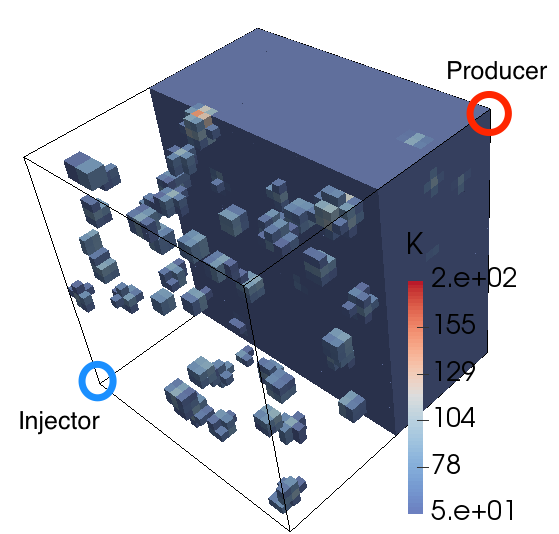}
\caption{Example 7. Setup with a random absolute permeabilities in a three dimensional domain. }
\label{fig:ex_8_set}
\end{figure}
In this example, we simply extend the previous example to a three dimensional domain $\Omega = (\SI{0},\SI{1}{\metre})^3$ with absolute permeabilities given as figure 
\ref{fig:ex_8_set}.
Wells are specified at the corners with 
injection at $(0,0,0)$ and production at $(\SI{1}{\metre},\SI{1}{\metre},\SI{1}{\metre})$. 
The numerical parameters are chosen as 
$h_{\min}=\num{5.4e-02}$ and $\Delta t = \num{3.4e-03}$.
All the other physical parameters and boundary conditions are the same as in the previous example.

\begin{figure}[!h]
\centering
\subfloat[$k=110$ ]
{
\includegraphics[width=0.3\textwidth]{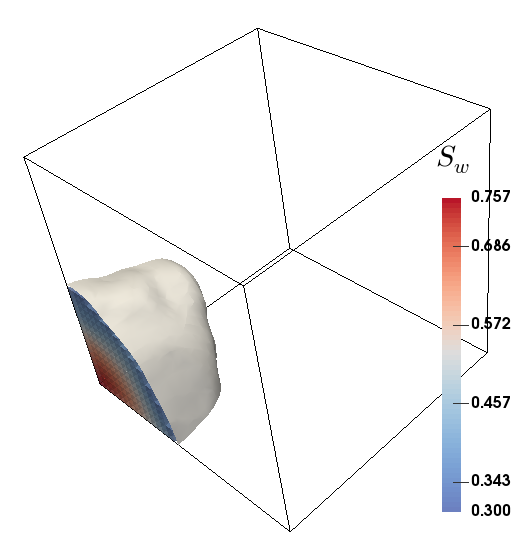}
}
\subfloat[$k=210$]
{
\includegraphics[width=0.3\textwidth]{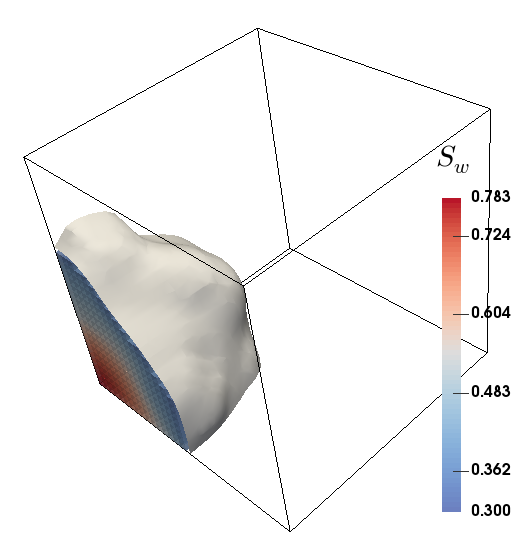}
}
\subfloat[$k=360$ ]
{
\includegraphics[width=0.3\textwidth]{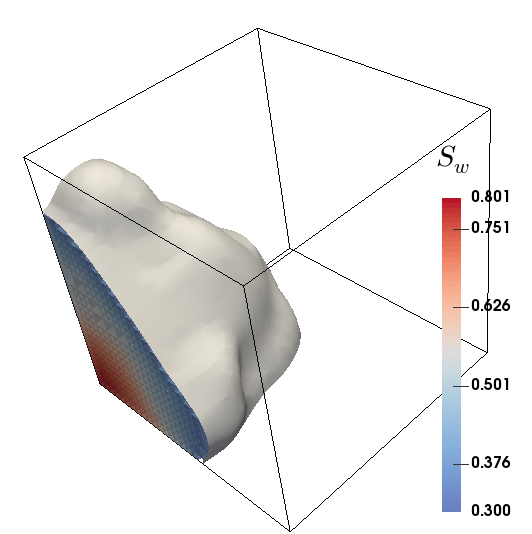}
}\\
\subfloat[$k=510$ ]
{
\includegraphics[width=0.3\textwidth]{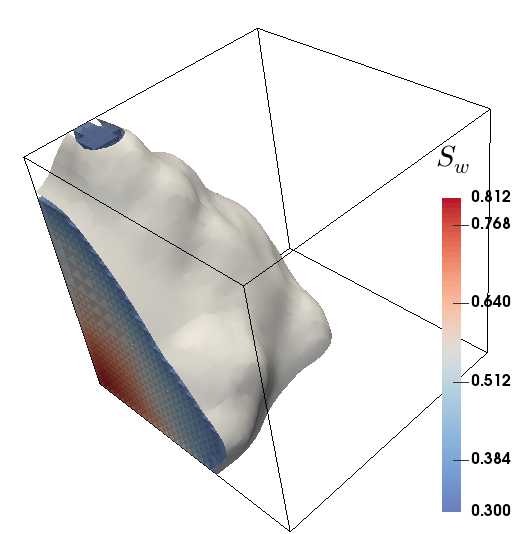}
}
\subfloat[$k=1010$ ]
{
\includegraphics[width=0.3\textwidth]{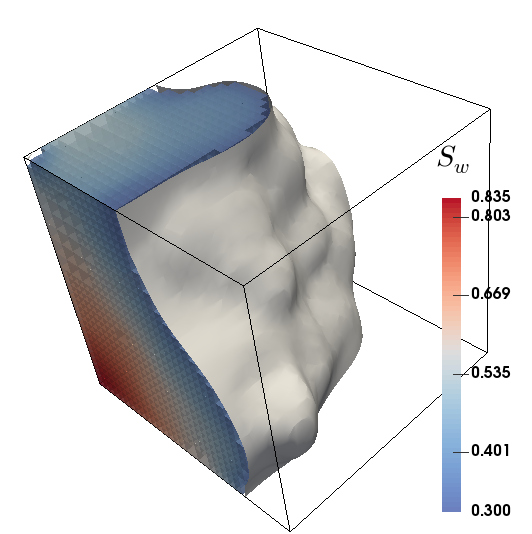}
}
\subfloat[$k=1290$ ]
{
\includegraphics[width=0.3\textwidth]{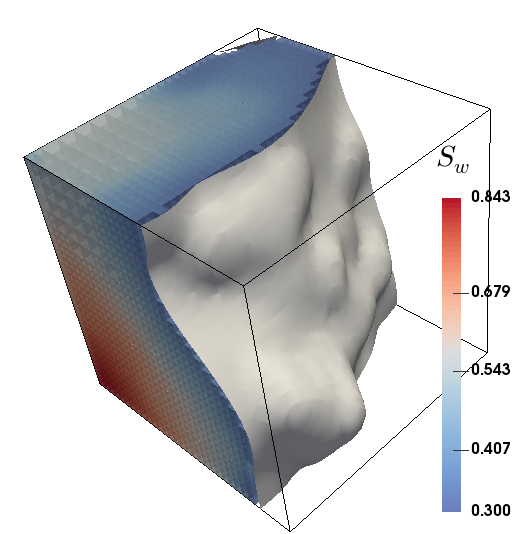}
}
\caption{Example 7. Contour value $S_w = 0.3$ for each time step. 
}
\label{fig:ex_7_3d}
\end{figure}

Figure \ref{fig:ex_7_3d} illustrates the contour value  of $S_w = 0.3$ for each time step. 
Here the maximum EG-$\polQ_1$ degrees of freedom for wetting phase saturation at the final time step is around $70,000$ and this example is computed by employing four multiple parallel processors (MPI).  

\subsection{Example 8. Well injections with gravity and a capillary pressure} 
\begin{figure}[!h]
\centering
\subfloat[Setup]
{
\includegraphics[width=0.3\textwidth]{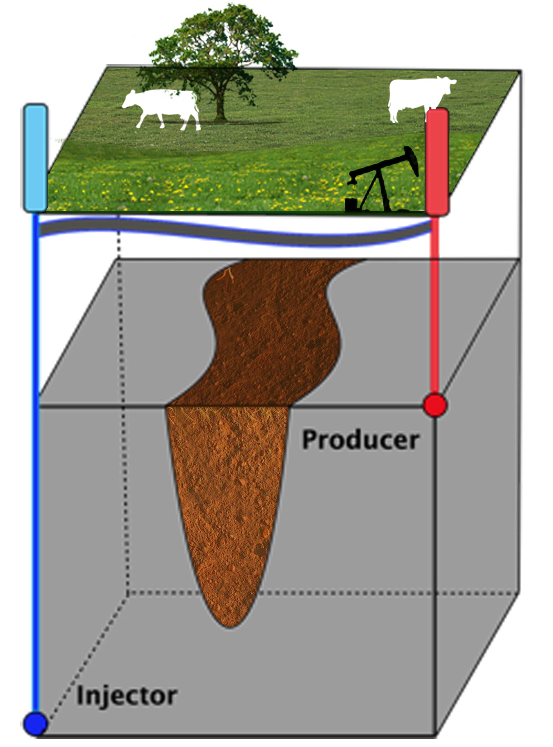}
\label{fig:ex4_a}
}
\hspace*{0.05in}
\subfloat[A domain]
{
\includegraphics[width=0.3\textwidth]{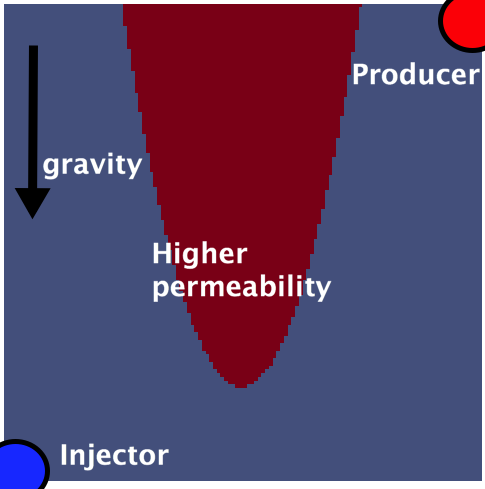}
\label{fig:ex4_b}
}
\caption{Example 6. Setup with the domain and the boundary conditions. (b) Two dimensional computational domain is defined by slicing the three dimensional domain (a) vertically. 
Bottom blue is the injection well and top red is the production well in the reservoir. Higher permeability zone is in the middle due to long sediments. }
\label{fig:ex_5_set}
\end{figure}
\begin{figure}[!h]
\centering
\subfloat[t=145s]
{
\includegraphics[width=0.25\textwidth]{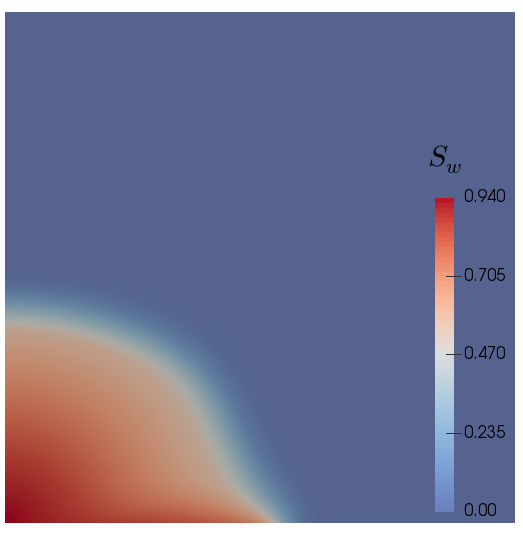}
}
\subfloat[t=250s]
{
\includegraphics[width=0.25\textwidth]{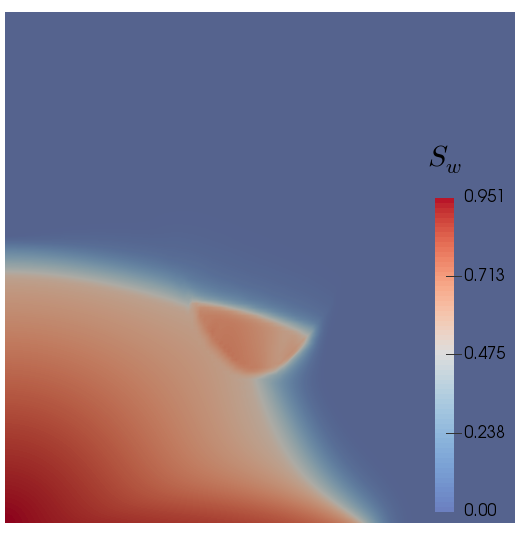}
}
\subfloat[t=325s]
{
\includegraphics[width=0.25\textwidth]{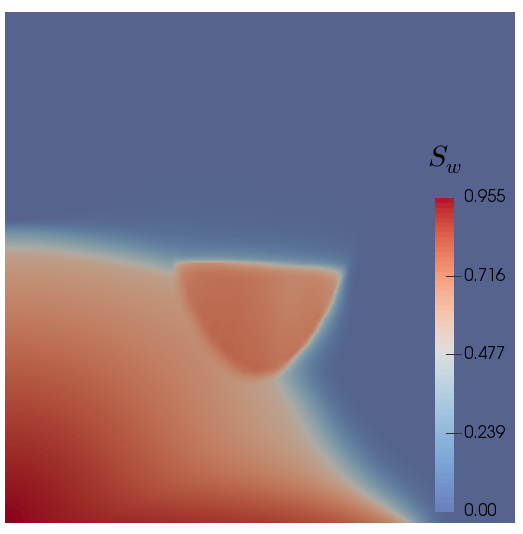}
} \\
\subfloat[t=500s]
{
\includegraphics[width=0.25\textwidth]{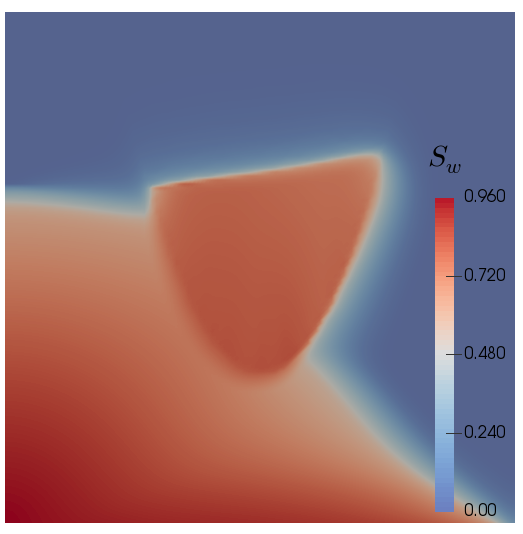}
}
\subfloat[t=685s]
{
\includegraphics[width=0.25\textwidth]{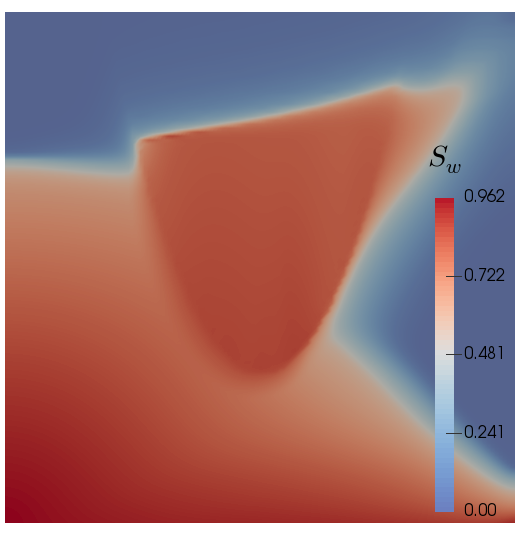}
}
\subfloat[t=720s]
{
\includegraphics[width=0.25\textwidth]{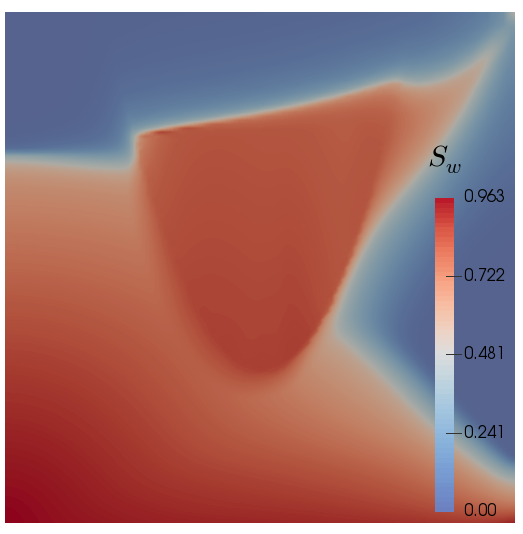}
}
\caption{Example 6. The wetting phase saturation ($S_w$) at each time.}
\label{fig:ex_5_num}
\end{figure}

Figure \ref{fig:ex4_a} illustrates an example of an existing reservoir 
where we have sliced a computational domain vertically,  
{$\Omega = (\SI{0}{\metre},\SI{50}{\metre})^2$} as shown in Figure \ref{fig:ex4_b}.
Wells are rate specified at the corners with 
injection at $(0,0)$ and production at $(\SI{50}{\metre},\SI{50}{\metre})$. 
A high permeability zone representing long sediments is located at ($y \geq 0.16x^2 - 7.78x + 112.22$), where $K_D=\SI{10}{D}$ and $K_D=\SI{1}{D}$ otherwise.
We assume the domain is saturated with a non-wetting phase, 
i.e $s^0_n = 1$ and $s^0_w = 0$ and  a wetting phase fluid  is injected.
Fluid and rock properties are given as 
$\mu_w = \SI{1}{cP}$, $\mu_n = \SI{3}{cP}$, 
$\rho_w = \SI{1000}{kg/m^3}$, $\rho_n=\SI{830}{kg/m^3}$, $c_w^F = \num{e-10}$, 
$f_w^+ = \SI{2.5}{m/s}$,  
$f_w^- = -\SI{2.5}{m/s}$,
$f_n = 0$,
and $\phi = 0.2$.
Relative permeabilities are given as functions of the wetting phase saturation \eqref{ex1_relperm}, and the capillary pressure is set with {$B_c = -0.001$ and $\varepsilon_s=0.1$}.
The penalty coefficients are set as $\alpha = 1$, $\alpha_c = 1$ and $\alpha_T = 1000$ and the time step is set by $\Delta t= \num{0.18}$.
Here, we employ the gravity $\bg =[0,\SI{-9.8}{m/s^2}$], and for the same scaling with pressure (atm), we divide it by ${101325}$ ($\SI{1}{atm} = \SI{101325}{\pascal}$).
Figure \ref{fig:ex_5_num} illustrates the injected wetting phase saturation values for each time step number. 
We observe the effect of the gravity.

\begin{figure}[!h]
\centering
\subfloat[$t=\SI{14.5}{s}$]
{
\includegraphics[width=0.25\textwidth]{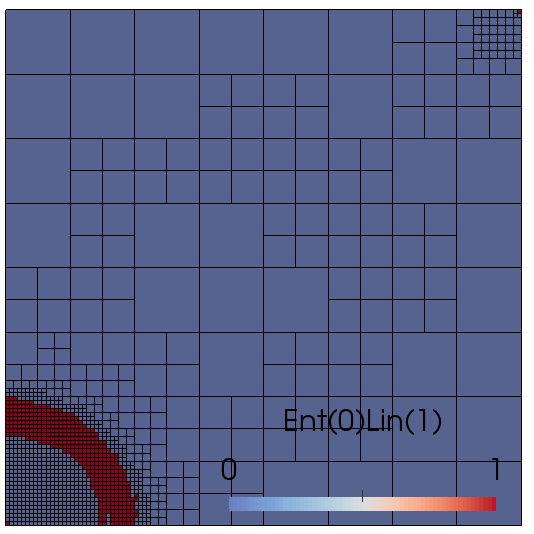}
}
\subfloat[$t=\SI{650}{s}$]
{
\includegraphics[width=0.25\textwidth]{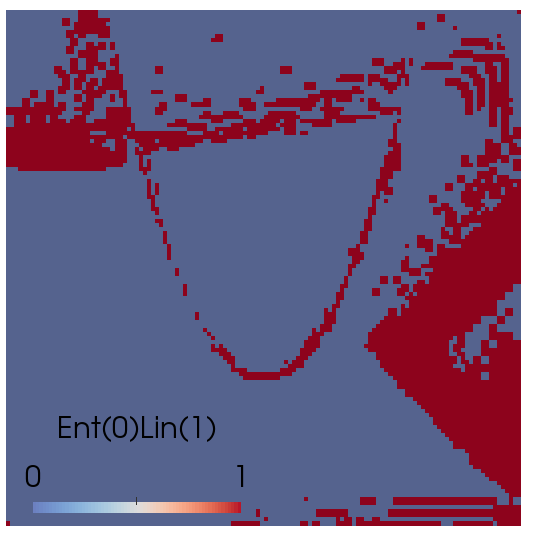}
}
\caption{Example 6. Entropy choices for early and later time.}
\label{fig:ex_5_num_ent}
\end{figure}

\begin{figure}[!h]
\centering
\subfloat[$S_n$ values in time]
{
\includegraphics[width=0.45\textwidth]{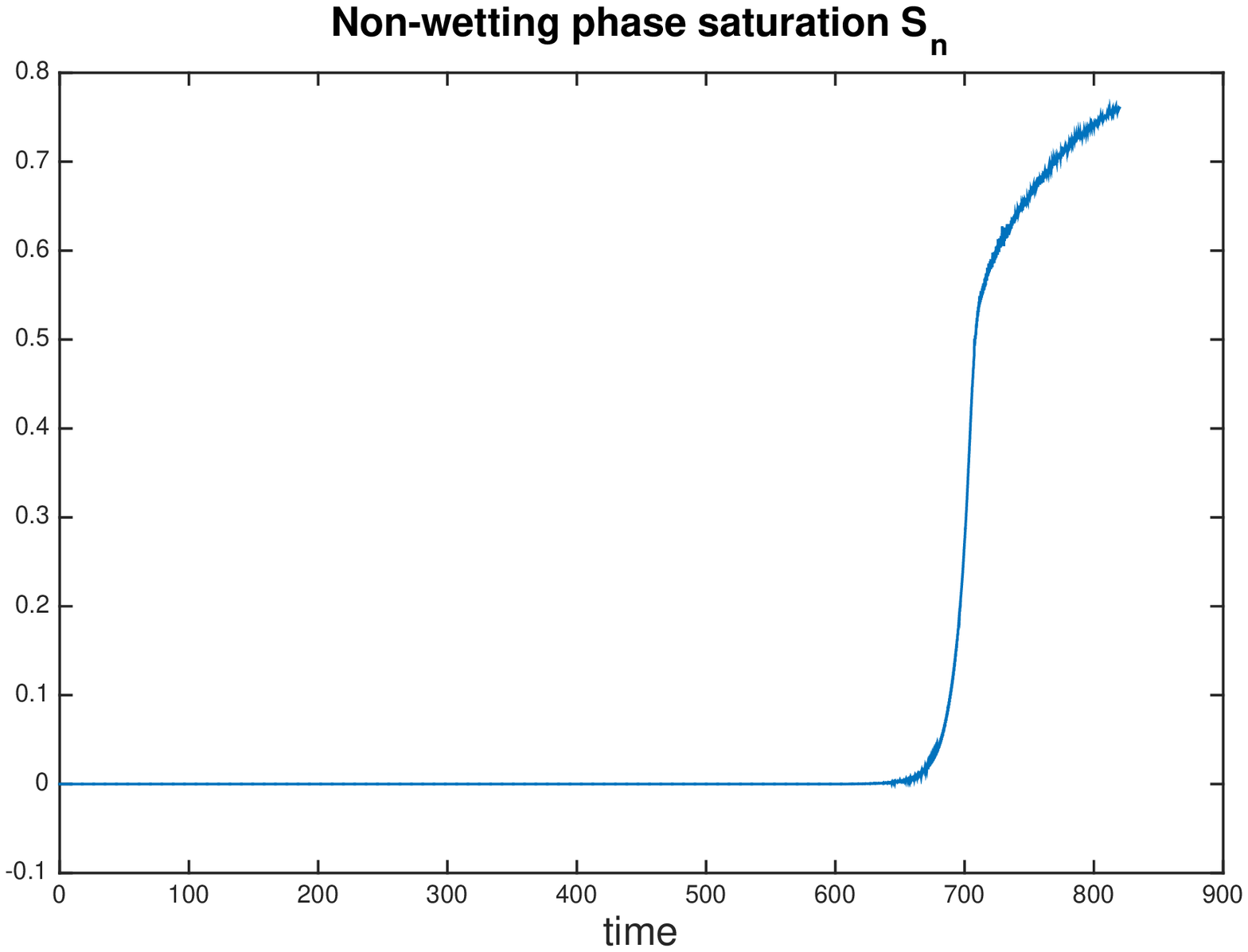}
\label{fig:ex_5_num_prod_a}
}
\subfloat[Accumulative oil production rate ($\sum f \times S_n$)] 
{
\includegraphics[width=0.45\textwidth]{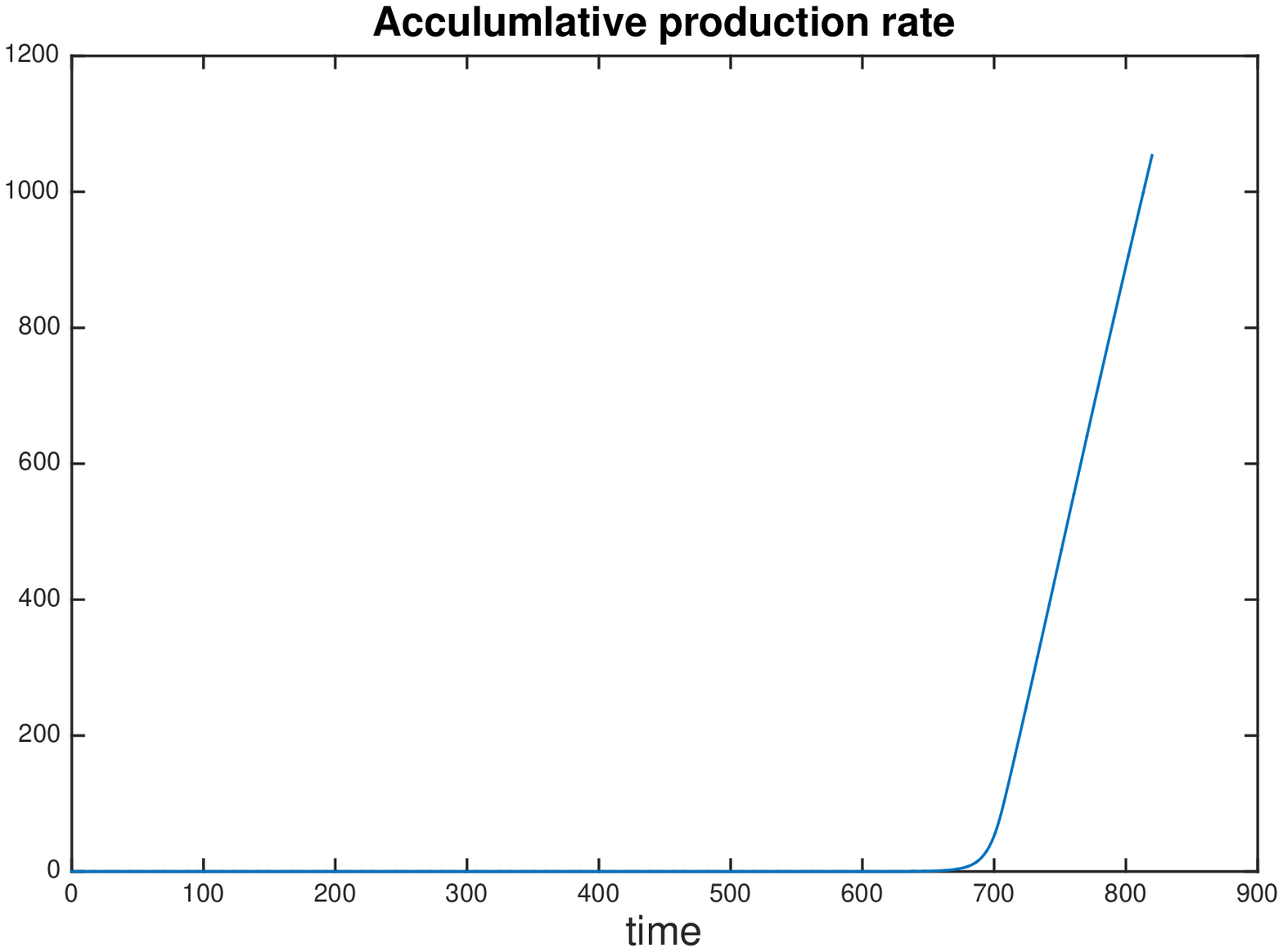}
\label{fig:ex_5_num_prod_b}
}
\caption{Example 6. Production data}
\label{fig:ex_5_prod}
\end{figure}

The entropy stabilization coefficients  are set as 
$\lambda_{\textsf{Ent}} = 40$ and 
$\lambda_{\textsf{Lin}}=1$,  
where the entropy function \eqref{eqn:entropy_func_2} is 
chosen with $\varepsilon = \num{e-3}$.
Figure \ref{fig:ex_5_num_ent} illustrates the choice for stabilization. 
Dynamic mesh adaptivity is employed with initial refinement level 
$\textsf{Ref}_T =4$,  $R_{\max}=7$ and  $R_{\min}=3$ with 
a minimum mesh size is $h_{\min}=\num{0.4}$.
In addition, Figure \ref{fig:ex_5_prod} presents the production data.  
The oil saturation values (non-wetting phase $S_n$) over the time are plotted with the accumulative oil production rate ($\sum_{k=0}^\mathbb{T} |S_n {f}^-|$).

\section{Conclusion}
\label{sec:conclusion}
In this paper, we present enriched Galerkin (EG) approximations for two-phase flow problems in porous media with capillary pressure. 
EG preserves local and global conservation for fluxes and has 
fewer degrees of freedom compared to DG. 
For a high order EG transport system, entropy residual stabilization is applied to avoid spurious oscillations. 
In addition, dynamic mesh adaptivity employing entropy residual as an error indicator reduces computational costs for large-scale computations.
Several examples in two and three dimensions including
error convergences and 
a well known capillary pressure benchmark problem are shown in order to verify and demonstrate the performance of the algorithm.
Additional challenging effects arising from gravity and rough relative permeabilities for foam are presented.

\section*{Acknowledgments}
The research by S. Lee and  M. F. Wheeler was partially supported by 
a DOE grant DE-FG02-04ER25617 and
Center for Frontiers of Subsurface Energy Security, an Energy Frontier
Research Center funded by the U.S. Department of Energy, Office of Science,
and Office of Basic Energy Sciences, DOE Project $\#$DE-SC0001114.
M. F. Wheeler was also partially supported by Moncrief Grand Challenge Faculty Awards from 
The Institute for Computational Engineering and Sciences (ICES), the University of Texas at Austin.

\bibliographystyle{spmpsci} 
\bibliography{xyz}

\end{document}